\documentclass[11pt]{amsart}

\newif\iffinalrun
\finalruntrue

\newif\ifkuvio

\iffinalrun
  \kuviotrue
\fi

\usepackage{amssymb,eucal}
\iffinalrun
\else
  \usepackage[notref,notcite]{showkeys}
\fi

\ifkuvio
  \usepackage[forcekdg,arrsy]{kuvio}
\fi

\usepackage{hyperref}

\newcommand{\tens}{\otimes}

\newcommand{\A}{\mathbb A}

\newcommand{\B}{\mathbb B}
\newcommand{\C}{\mathbb C}

\newcommand{\G}{\mathbb G}
\newcommand{\tG}{\widetilde G}
\newcommand{\HH}{\mathcal H}
\newcommand{\tHH}{\widetilde\HH}
\renewcommand{\k}{{\bar k}}
\renewcommand{\L}{\mathbb L}
\newcommand{\tL}{\widetilde L}

\newcommand{\N}{\mathbb N}
\newcommand{\oN}{\overline N}
\newcommand{\tN}{\widetilde N}
\renewcommand{\O}{{\mathcal O}}

\newcommand{\oP}{\overline P}
\newcommand{\Q}{\mathbb Q}

\newcommand{\qp}{\Q_p}

\newcommand{\R}{\mathbb R}
\renewcommand{\SS}{\mathcal S}
\newcommand{\tSS}{\widetilde\SS}
\newcommand{\tS}{\widetilde S}
\newcommand{\T}{\mathbb T}
\newcommand{\tT}{\widetilde T}
\newcommand{\ttau}{\widetilde\tau}
\newcommand{\oU}{{\overline U}}

\newcommand{\oW}{{\overline W}}
\newcommand{\hS}{\widehat S}
\newcommand{\Z}{\mathbb Z}

\newcommand{\oPhi}{\overline\Phi}

\newcommand{\INTO}{\hookrightarrow}
\newcommand{\onto}{\twoheadrightarrow}
\newcommand{\congto}{\xrightarrow{\,\sim\,}}
\newcommand{\tocong}{\xleftarrow{\,\sim\,}}

\newcommand{\s}{^\times}

\newcommand{\dual}{^\vee}

\DeclareMathOperator{\Aut}{Aut}
\DeclareMathOperator{\End}{End}
\DeclareMathOperator{\Hom}{Hom}
\DeclareMathOperator{\GL}{GL}
\DeclareMathOperator{\Spec}{Spec}

\newcommand{\Gal}{\mathop{\mathrm{Gal}}}
\newcommand{\Stab}{\mathop{\mathrm{Stab}}}
\newcommand{\tr}{\operatorname{tr}}

\newcommand{\ord}{\operatorname{ord}}

\newcommand{\nr}{^\mathrm{nr}}
\newcommand{\nd}{_\mathrm{nd}}
\newcommand{\Aff}{\mathop{\mathrm{Aff}}}

\DeclareMathOperator{\Lie}{Lie}
\DeclareMathOperator{\rank}{rank}
\DeclareMathOperator{\diag}{diag}

\newcommand{\ind}{\operatorname{c-Ind}}

\newcommand{\red}{{\operatorname{red}}}

\newcommand{\im}{\operatorname{im}}

\newcommand{\supp}{\operatorname{supp}}

\iffinalrun
  \newcommand{\need}[1]{}
  \newcommand{\mar}[1]{}
\else
  \newcommand{\need}[1]{{\tiny *** #1}}
  \newcommand{\mar}[1]{\marginpar{\tiny #1}}
\fi

\newcommand{\rr}{\rrbracket}

\renewcommand{\)}{\textup{)}}

\theoremstyle{plain} 
\newtheorem{lm}[equation]{Lemma}

\newtheorem{prop}[equation]{Proposition}
\newtheorem{thm}[equation]{Theorem}
\newtheorem{coroll}[equation]{Corollary}

\theoremstyle{definition}
\newtheorem{df}[equation]{Definition}
\newtheorem{rk}[equation]{Remark}

\numberwithin{equation}{section}
\numberwithin{figure}{section}

\begin{document}

\title[A Satake isomorphism in characteristic $p$]{A Satake isomorphism in characteristic $p$}
\author{Florian Herzig}
\address{Department of Mathematics, 2033 Sheridan Road, Evanston IL 60208-2730, USA}
\email{herzig@math.northwestern.edu}
\thanks{The author was partially supported by NSF grant DMS-0902044.}
\maketitle

\begin{abstract}
  Suppose that~$G$ is a connected reductive group over a $p$-adic field~$F$, that $K$ is a hyperspecial maximal compact
  subgroup of~$G(F)$, and that $V$ is an irreducible representation of~$K$ over the algebraic closure of the residue field of~$F$.
  We establish an analogue of the Satake isomorphism for the Hecke algebra of compactly supported, $K$-biequivariant
  functions $f: G(F) \to \End V$. These Hecke algebras were first considered by Barthel--Livn\'e for~$\GL_2$. They
  play a role in the recent mod~$p$ and $p$-adic Langlands correspondences for $\GL_2(\qp)$, in generalisations
  of Serre's conjecture on the modularity of mod~$p$ Galois representations, and in the classification of irreducible mod~$p$
  representations of unramified $p$-adic reductive groups.
\end{abstract}

\section{Introduction}

\subsection{Statement of the theorem}

Let~$F$ be a finite extension of~$\qp$ with ring of integers~$\O$, uniformiser~$\varpi$, and residue field~$k$ of order~$q$.
Suppose that $G$ is a connected reductive group over~$F$ that is unramified (i.e., quasi-split and split over an unramified extension)
and that $K$ is a hyperspecial maximal compact subgroup. Fix any maximal split torus~$S$ in~$G$ 
such that the apartment corresponding to~$S$ contains the hyperspecial point in the reduced building corresponding to~$K$. Since $G$
is quasi-split, $T = Z_G(S)$ is a maximal torus of~$G$.

With these assumptions it is known that $G$ extends to a smooth $\O$-group scheme \cite[\S 3.8.1]{bib:Tits}, that we will also denote by~$G$, whose
special fibre is a connected reductive group over~$k$ and such that $K = G(\O)$. The tori $S$ and~$T$ extend to smooth
$\O$-subgroup schemes $S \subset T$ of~$G$ which reduce to a maximal split torus and its centraliser in the special fibre
of~$G$.  The relative root systems of~$S$ in~$G$ in the two fibres are naturally identified with each other. We denote by
$\Phi \subset X^*(S)$ the set of roots, by~$\Phi^+$ a choice of positive roots, and by~$W$ the Weyl group.  There is a closed
$\O$-subgroup scheme $B = T \ltimes U$ of~$G$ whose fibres are the Borel subgroups associated to~$\Phi^+$.

Suppose that $V$ is an irreducible representation of~$G(k)$ over~$\k$, which we also consider as a representation of~$K$ via
the reduction homomorphism $K = G(\O) \to G(k)$. The \emph{Hecke algebra $\HH_G(V)$ of~$V$} is the $\k$-algebra of compactly
supported functions $f : G(F)\to \End_{\k} V$ satisfying $f(k_1gk_2) = k_1 f(g)k_2$ for all $k_1$, $k_2 \in K$, $g \in G(F)$,
where the multiplication is given by convolution. We remark that by Frobenius reciprocity it follows that $\HH_G(V) \cong
\End_{\k G(F)} (\ind_K^{G(F)} V)$, where $\ind_K^{G(F)} V = \{ \psi : G(F) \to V : \psi(kg) = k\psi(g)\; \forall k \in K\;
\forall g \in G;\; \text{$\supp \psi$ compact} \}$ is the compactly induced representation (see \cite[Prop. 5]{bib:BL2}).

It is known that the $T(k)$-representation $V^{U(k)}$ is one-dimensional (see Lemma~\ref{lm:U-inv}). The corresponding Hecke
algebra $\HH_T(V^{U(k)})$ consists of $T(\O)$-biequivariant, compactly supported functions $\varphi: T(F) \to \End_{\k}(V^{U(k)}) = \k$.

Let $\ord_F : F\s \onto \Z$ denote the valuation of $F$.
For $\chi \in X^*(S)$ and $t \in T(F)$ we define $(\ord_F \circ \chi)(t)$ to be $\frac 1n\ord_F(n\chi(t))$, where
$n > 0$ is chosen so that $n\chi$~extends to an $F$-rational character of~$T$. Since $X^*_F(T) \to X^*(S)$ is injective
with finite cokernel, this does not depend on any choices.

\begin{df}\label{df:t-minus}
  Let~$T^-$ denote the following submonoid of~$T(F)$:
  \begin{equation*}
    T^- = \{t \in T(F) : (\ord_F\circ \alpha)(t) \le 0 \quad \forall \alpha \in \Phi^+\}.
  \end{equation*}
  Let $\HH^-_T(V^{U(k)})$ denote the subalgebra of~$\HH_T(V^{U(k)})$ consisting of those $\varphi : T(F) \to \k$ that are supported on~$T^-$.
\end{df}

\begin{thm}\label{thm:satake} Suppose that $V$ is an irreducible representation of~$G(k)$ over~$\k$.
  Then
  \begin{align*}
    \SS : \HH_G(V) &\to \HH_T(V^{U(k)})\\ 
    f &\mapsto \left(t \mapsto \sum_{u \in U(F)/U(\O)} f(tu)\Big|_{V^{U(k)}}\right)
  \end{align*}
  is an injective $\k$-algebra homomorphism with image $\HH^-_T(V^{U(k)})$.
\end{thm}

Note that as $f$ is compactly supported, the sum over $U(F)/U(\O)$ has only finitely many non-zero terms and $\SS f$ is compactly
supported. Since $T(F)$ normalises~$U(F)$, the image of~$\SS f$ is contained in~$V^{U(k)}$.

It is easy to see that $\lambda \mapsto \lambda(\varpi)$ yields an isomorphism $X_*(S) \to
T(F)/T(\O)$ which sends the antidominant coweights
$X_*(S)_- = \{ \lambda \in X_*(S) : \langle \lambda, \alpha \rangle \le 0\ \forall \alpha \in \Phi^+ \}$ 
to $T^-/T(\O)$ (Lemma~\ref{lm:torus-cosets}).

\begin{coroll}\label{cor:hecke_struct}
  $\HH_G(V)$ is commutative and isomorphic to $\k[X_*(S)_-]$. In particular it is noetherian.
\end{coroll}

At least when $G$ is split and the derived subgroup of~$G$ is simply connected, there is another argument to see that $\HH_G(V)$ is
commutative which uses an analogue of a Gelfand involution.  See the end of Section~\ref{sec:satake}.

\subsection{Comparison with the classical Satake isomorphism}

Recall that the classical Satake isomorphism is given by the formula
\begin{align*}
  \C[K\backslash G(F)/K] &\congto \C[X_*(S)]^W \\
  f &\mapsto \left(t \mapsto \delta(t)^{1/2} \int_{U(F)} f(tu) du \right),
\end{align*}
where $\delta$ is the modulus character of the Borel and the Haar measure~$du$ on~$U(F)$ satisfies $\int_{U(\O)} du = 1$
(\cite{bib:Cartier}, \cite{bib:Gross_Satake}). The relevance of the factor $\delta^{1/2}$ is to make the image of the Satake
transform $W$-invariant. Leaving it out still yields an algebra homomorphism~$\SS'$ into $\C[X_*(S)]$, which now also makes
sense over~$\Z$, and which is obviously compatible with~$\SS$ when $V$ is the trivial representation:
\ifkuvio
  \[ \cellwidth=45pt \Diag
  &\Z[K\backslash G(F)/K] & \rInto^{\SS'} & \Z[X_*(S)] \\
  &\dTo  && \dTo  \dy{-3mm} \\
  \HH_G(1) = {} & \k[K\backslash G(F)/K] & \rInto^{\SS} & \k[X_*(S)] \\ 
  \endDiag \]
\fi
In this case (when $V$ is trivial) there is a simple explanation why the image of~$\SS$ is supported on antidominant coweights.
The image of the Satake transform is $W$-invariant and the modulus character is a power of~$p$ which, among the $W$-conjugates
of a given coweight, is biggest on the antidominant one.

The proof of Theorem~\ref{thm:satake} follows the same steps as the classical proof, but there are two complications.
Firstly, it is harder to determine the space of Hecke operators supported on a given double coset. This requires an argument
using the Bruhat--Tits building (Prop.~\ref{prop:building}). Secondly, for general~$V$ it is subtle to prove that the image
of~$\SS$ is contained in $\HH_T^-(V^{U(k)})$. We first show that the image is supported on ``almost antidominant'' coweights
and then use that $\SS$ is a homomorphism to conclude. This extra step is really necessary, as one sees by considering
the Hecke bimodule $\Hom_{G(F)}(\ind_K^{G(F)} V_1,\ind_K^{G(F)} V_2)$ whose support under the Satake map may
extend slightly beyond the antidominant coweights \cite[\S 6]{bib:parab}.

\subsection{Comparison with the $p$-adic Satake isomorphism}

Schneider--Teitelbaum \cite{bib:ST} constructed $p$-adic Satake maps, and their $p$-adic completions, for the Hecke algebras
associated to an irreducible representation of~$G_{/F}$. In Proposition~\ref{prop:st-comparison} we establish a compatibility
between their $p$-adic Satake map and the mod~$p$ Satake map~$\SS$, in case $V$ extends to a representation of~$G_{/k}$.
(This is satisfied, for example, if the derived subgroup of~$G_{/\k}$ is simply connected.)  In this case $V$ is
a submodule of the reduction of a $K$-stable lattice in some irreducible representation of~$G_{/F}$. Note that
$V$ does not necessarily equal the reduction; in fact, this cannot usually be achieved.

\subsection{The $W$-regular case}

The refined Cartan decomposition says that the $\lambda(\varpi)$ for $\lambda \in X_*(S)_-$ form
a system of coset representatives for $K\backslash G(F)/K$.  We will see in the proof of Theorem~\ref{thm:satake} that
$\HH_G(V)$ has a natural $\k$-basis $\{T_\lambda : \lambda \in X_*(S)_- \}$. The Hecke operator $T_\lambda$ is characterised by having support
$K\lambda(\varpi)K$ and by $T_\lambda(\lambda(\varpi)) \in \End_\k V$ being a projection. More obviously
(Lemma~\ref{lm:torus-cosets}), $\HH^-_T(V^{U(k)})$ has a $\k$-basis $\{\tau_\lambda : \lambda \in X_*(S)_-\}$ where
$\tau_\lambda$ is supported on $\lambda(\varpi)T(\O)$ and $\tau_\lambda(\lambda(\varpi)) = 1$.

We will say that an irreducible representation $V$ of~$G(k)$ over~$\bar
k$ is \emph{$W$-regular} if the ``extremal weight subspaces'' $w V^{U(k)} \subset V$ for $w \in W$ are distinct.

\begin{prop}\label{prop:satake_regular}
  Suppose that $V$ is $W$-regular. Then for each $\lambda \in X_*(S)_-$, $\SS T_\lambda = \tau_\lambda$.
  In particular, $T_\lambda * T_{\lambda'} =
  T_{\lambda + \lambda'}$ for all $\lambda$, $\lambda' \in X_*(S)_-$.
\end{prop}

For general~$V$ and for $\lambda \in X_*(S)_-$, the proof of Theorem~\ref{thm:satake} shows that 
\begin{equation*}
  \tau_\lambda = \sum_{{\mu \in X_*(S)_-} \atop {\mu \ge_\R \lambda}} d_\lambda(\mu) \SS T_\mu,
\end{equation*}
where $d_\lambda(\mu) \in \k$ and $d_\lambda(\lambda) = 1$. In the classical setting the work of Lusztig and Kato shows that
the $d_\lambda(\mu)$ are Kazhdan--Lusztig polynomials in $q = |k|$. (See \cite{bib:Kato}, \cite{bib:Gross_Satake},
\cite{bib:HKP}.)  In~\cite[\S 5]{bib:parab} we use their results to compute $d_\lambda(\mu)$ in all cases,
at least when $G$ is split and has simply connected derived subgroup. It turns out that $d_\lambda(\mu)$ does not depend
on~$V$ but only on the stabiliser of the subspace $V^{U(k)}$ in~$W$.

\subsection{Satake parameters}

Let $\Delta \subset \Phi^+$
denote the set of simple roots. Let~$\hS$ be the torus dual to~$S$ (over~$\k$). For each
subset $J \subset \Delta$, define the torus $\hS_J$ by the exact sequence
\begin{equation*}
  \G_m^J \to \hS \to \hS_J \to 1,
\end{equation*}
where the first map is given by $\prod_{\delta \in J} \delta$. 
The closed points of the ``toric'' variety $\Spec \HH_G(V)$ have the following concrete description.
Classically only one torus ($\hS = \hS_\varnothing$) is needed.

\begin{coroll}\label{cor:satake_params}
  The $\k$-algebra homomorphisms $\HH_G(V) \to \k$ are parameterised by pairs $(J,s_J)$, where $J \subset \Delta$ and
$s_J \in \hS_J(\k)$.
\end{coroll}

\mar{It's then easy to describe the compatibility with changing the uniformiser or the choice of Borel. To be added?}

In~\cite[\S 4]{bib:parab} we give an alternative parameterisation, analogous to the classical parameterisation by unramified
characters of~$T$.

\subsection{Example: $G = \GL_n$} We suppose that $S = T$ is the diagonal torus and that $B$ is the Borel subgroup of
upper-triangular matrices.  Then the $\lambda_i(x) = \diag(1,\dots,1,x,\dots,x)$ (with $i$ non-trivial entries) generate
$X_*(S)_-$, and we denote by~$T_i$ the corresponding Hecke operator $T_{\lambda_i}$. Theorem~\ref{thm:satake} shows that
$\HH_G(V)$ is the localised polynomial algebra $\k[T_1,\dots,T_{n-1},T_n^{\pm 1}]$. 

\subsection{Previous work} The Hecke algebras $\HH_G(V)$ were first calculated by Barthel--Livn\'e when $G = \GL_2$ \cite{bib:BL1},
\cite{bib:BL2}. (We follow their strategy for computing $\HH_G(V)$ as vector space. However they used explicit methods to
determine the algebra structure.)  This was important for their (partial) classification of irreducible smooth
representations~$\pi$ of~$\GL_2(F)$ over~$\k$ that have a central character, which was completed by Breuil when $F=\qp$
\cite{bib:Breuil} and which plays a crucial role for mod~$p$ and $p$-adic local Langlands correspondences for $\GL_2(\qp)$.
In~\cite{bib:parab} we extend the work of Barthel--Livn\'e, giving a classification of irreducible admissible
representations of~$\GL_n(F)$ over $\k$, in terms of supersingular representations.  Our proofs heavily depend on the methods
developed in this paper.

We also remark that Schein independently determined the Hecke algebras for~$\GL_n$ by explicit methods \cite{bib:Schein},
after we had done this in a similar manner.

In another direction, Gross showed that the classical Satake isomorphism can be defined over~$\Z[q^{1/2}, q^{-1/2}]$
\cite[\S 3]{bib:Gross_Satake}. See also \cite[\S 1.2]{bib:lazarus}.

\subsection{Algebraic modular forms}

Suppose that $F = \qp$ and that $G$ arises by base change from a connected reductive $\Q$-group~$G$ such that $G(\R)$ is
compact. Given a compact open subgroup $K_{\A} = K \times K^p$ in~$G(\A^\infty)$ we can consider Gross's space $M(K_\A,V^*)$
of algebraic modular forms of level $K_{\A}$ and weight~$V^*$, the linear dual of~$V$ \cite{bib:Gross_AMF}. The Hecke
algebra~$\HH_G(V)$ naturally acts on it and there is a simple compatibility result of the action of the~$T_\lambda$ on
$M(K_\A,V^*)$ with classical Hecke operators.
In joint work with Matthew Emerton and Toby Gee we use it to prove strong new results on the weights in a Serre-type conjecture
for rank~3 unitary groups
\cite{bib:EGH}.

\subsection{Organisation of the paper}

In Section~\ref{sec:proofs} we discuss the proofs of the main results.
Technical parts of the arguments requiring buildings are discussed in
Section~\ref{sec:buildings}. We include the proofs of some well-known results since we did not find an appropriate reference.

For a reader who is inexperienced with algebraic groups, we recommend to assume first that $G = \GL_n$ or, more generally, that
$G$ is split with simply connected derived subgroup. Many arguments simplify in these settings.

\subsection{Acknowledgements} I thank Kevin Buzzard for suggesting the problem of computing these Hecke algebras in the case
of $G = \GL_3$. I thank Matthew Emerton for encouraging me to generalise and for many helpful discussions. I am grateful to
Jiu-Kang Yu for his helpful comments related to Section~\ref{sec:buildings} and to Vytautas Pa\v sk\= unas and Peter Schneider for
helpful discussions. I thank the referee for a careful reading of the paper.

\section{Proofs}\label{sec:proofs}

\subsection{The Satake isomorphism for $\HH_G(V)$}\label{sec:satake}

\begin{lm}\label{lm:torus-cosets}
  The map $\zeta : T(F) \to X_*(T)$ given by
  \begin{equation*}
    \langle \zeta(t),\chi \rangle = \ord_F(\chi(t)) \quad \forall \chi \in X^*(T)
  \end{equation*}
  induces isomorphisms of abelian groups
  \begin{equation}\label{eq:6}
    S(F)/S(\O) \congto T(F)/T(\O) \congto X_*(S).
  \end{equation}
  Moreover $T^-/T(\O)$ \(see Def.~\ref{df:t-minus}\) corresponds to $X_*(S)_-$ under the isomorphism.
  A ``splitting'' of~\eqref{eq:6} is provided by $X_*(S) \to S(F)$, $\lambda \mapsto \lambda(\varpi)$.
\end{lm}

Note that $\chi(t) \in (F\nr)\s$ since $T$ splits over an unramified extension, so $\ord_F(\chi(t)) \in \Z$.

\begin{proof}
  We consider the following diagram.
  \ifkuvio
  \[ \cellwidth=45pt \Diag
  0 & \rTo & S(\O) & \rTo & S(F) & \rTo & X_*(S) & \rTo & 0 \\
  && \dInto && \dInto  && \dInto \dy{-3mm} \\
  0 & \rTo & T(\O) & \rTo & T(F) & \rTo^{\zeta} & X_*(T) \\
  \endDiag \]
  \fi
  Note that $\zeta$ lands in the $\Gal(\overline F/F)$-invariant part of~$X_*(T)$, that is in~$X_*(S)$.
  As $T_{/\O\nr}$ is split (see Lemma~\ref{lm:conn_torus}), $\ker \zeta = T(F) \cap T(\O\nr) = T(\O)$. All the claims follow immediately.
\end{proof}

We need to introduce a partial order $\le_\R$ on $X_*(S)_\R$. First note that $X^*(S)_\R = \R\langle \Phi\rangle \oplus
X^*(G_{/F})_\R$,\mar{more details??} where $X^*(G_{/F}) = \Hom_F(G_{/F},\G_m)$. Since $\Phi$ is a root system in $\R\langle \Phi\rangle$, for every $\alpha
\in \Phi$ there is a ``coroot'' $\alpha\dual \in (\R\langle \Phi\rangle)^*$, characterised by $s_\alpha(x) = x-\langle x,\alpha\dual\rangle\alpha$.
We say that $y \ge_\R y'$ for $y$, $y' \in
(\R\langle \Phi\rangle)^*$ if $y-y'$ is a non-negative real linear combination of the positive coroots.

\begin{df}\label{df:partial-order}
  Suppose that $\lambda$, $\lambda' \in X_*(S)_\R$. We say that $\lambda \ge_\R \lambda'$ if
  $\lambda-\lambda'$ lies in the direct summand $(\R\langle \Phi\rangle)^*$ and $\lambda-\lambda' \ge_\R 0$.
\end{df}

Alternatively one could use the relative coroots in~$X_*(S)$ as defined in \cite[\S 15.3]{bib:Springer_LAG}.

\begin{lm}\label{lm:partial-order}
  Suppose that $\lambda \in X_*(S)$. Then $\{\lambda' \in X_*(S)_- : \lambda' \ge_\R \lambda \}$ is finite.
\end{lm}

\begin{proof}
  By the definition of~$\ge_\R$ we may project onto $(\R\langle \Phi\rangle)^*$. The projections $\bar\lambda$, $\bar\lambda'$ lie in the
  coweight lattice for the root system~$\Phi$ in $(\R\langle \Phi\rangle)^*$ and $\bar\lambda$ is antidominant.
  In this setting the result is well known.
\end{proof}

Next we will study the invariants of an irreducible $G(k)$-representation~$V$ over~$\k$ under the unipotent radical of a
parabolic subgroup.

\begin{lm}\label{lm:U-inv}
  Suppose that $V$ is an irreducible representation of~$G(k)$ over~$\k$. Then $V^{U(k)}$ is one-dimensional.
  Suppose that $P = L \ltimes N$ is a parabolic
  subgroup of~$G_{/k}$ and denote by $\oP = L\ltimes \oN$ the opposite parabolic.
  \begin{enumerate}
  \item $V^{N(k)}$ is an irreducible representation of~$L(k)$.
  \item The natural map $V^{N(k)} \to V \to V_{\oN(k)}$ is an isomorphism of $L(k)$-representations.
  \end{enumerate}
\end{lm}

Part~(i) was first proved by Smith \cite{bib:Smith} in the case that $G_{/k}$ is semisimple and simply connected.
Cabanes provided a general proof \cite{bib:Cabanes}, using $(B,N)$-pairs.
Below we give a proof that generalises the proof in the simply connected case found in~\cite[\S 5.10]{bib:Humph_LieType}.

\begin{proof}
  \emph{Let us first assume that the derived subgroup of~$G_{/\k}$ is simply connected.} By conjugating, we may assume that
  $P = L\ltimes N$ is a standard Levi decomposition, i.e., $P \supset B$ and $T \subset L$.  Let~$\G$ be the split $k$-form
  of~$G_{/\k}$ and fix a split maximal torus~$\T$ and a Borel~$\B$ containing it. Let $\phi \in \Gal(\k/k)$ denote the
  Frobenius element. There is a finite-order automorphism $\pi \in \Aut_k(\G,\B,\T)$ and an isomorphism $f : G(\k) \to
  \G(\k)$ respecting maximal tori and Borel subgroups such that $f \circ \phi = (\pi \circ \phi) \circ f$.  In particular,
  $G(k) = \G(\k)^{\pi \circ \phi}$.
  Let $\L \ltimes \N$ be the parabolic subgroup of~$\G$ corresponding to $L \ltimes N$ in~$G$.
  
  Since $\G'$ is simply connected, a (slight extension of a) result of Steinberg shows that $V$ is isomorphic to the
  restriction to~$G(k)$ of an irreducible representation $F(\nu)$ of the algebraic group~$\G$ whose highest weight $\nu \in
  X^*(\T)$ is $q$-restricted, i.e., satisfies $0 \le \langle \nu, \beta\dual\rangle < q$ for all simple
  roots~$\beta$ of~$\G$ \cite[Prop.~A.1.3]{bib:thesis}. Moreover $V^{U(k)} \cong F(\nu)_\nu$ 
  (the weight space of weight $\nu$) is one-dimensional.
  
  (i) This is Cor.~5.10 in~\cite{bib:Humph_LieType}. Even though $\G$ is assumed to be semisimple in that reference, the
  proof goes through word for word. From the proof we see that $F(\nu)^{\N} = F(\nu)^{N(k)}$ is the sum of weight spaces
  $F(\nu)_{\nu'}$ with $\nu-\nu' \in \Z_{\ge 0} \Theta^+$, where $\Theta^+$ are the positive roots of~$(\T,\L)$.  This is an
  irreducible $L(k)$-representation since $\nu$ is also $q$-restricted for~$\L$ and
  since $\L'$ is simply connected (as $\G'$ is simply connected). \mar{Actually by taking a
    $z$-extension as in the proof in the general case, we see that the irreducibility does not need the derived subgroup
    to be s.c.}

  (ii) Since $(V^*)^{\oN(k)} \cong \Hom_{\k}(V_{\oN(k)},\k)$ it follows that $V_{\oN(k)} \cong ((V^*)^{\oN(k)})^*$ is irreducible
  as $L(k)$-representation. It thus suffices to show that $V^{N(k)} \to V_{\oN(k)}$ is non-zero, or equivalently that
  $V^{N(k)}$ pairs non-trivially with $(V^*)^{\oN(k)}$ under the duality $V \times V^* \to \k$.
  By part~(i), $V^{N(k)}$ contains the highest weight space $L(\nu)_\nu$ and $(V^*)^{\oN(k)}$ contains the lowest weight space $(L(\nu)^*)_{-\nu}$.
  Since these pair non-trivially, this completes the proof. (One even sees directly in this way that the pairing on
  $V^{N(k)} \times (V^*)^{\oN(k)}$ is non-degenerate, i.e., that the map $V^{N(k)} \to V_{\oN(k)}$ is an isomorphism.)

  We remark that this argument shows that $V^{N(k)}$ is a direct summand of~$V$ as $L(k)$-representation, which is also clear from
  the proof in~\cite{bib:Humph_LieType}.
  
  \emph{Let us now reduce the general case to the previous one.} For ease of notation we will be writing~$G$ for its special
  fibre $G \times_\O k$, and similarly for $S$, $T$, etc.  We pick a $z$-extension of~$G$. This is an exact sequence
  \begin{equation}\label{eq:4}
    1 \to R \to \tG \xrightarrow{\pi} G \to 1
  \end{equation}
  of affine algebraic $k$-groups, where $\tG$ is reductive with $\tG'$ simply connected and $R$ a central torus (even an
  induced torus). Exactness means that the first map is a closed embedding, the second map faithfully flat, and that the
  first map is the kernel of the second. The notion of a $z$-extension goes back to Langlands in characteristic zero; for
  the general case see~\cite[\S 3.1]{bib:CT}.
  
  By~\cite[Thm.~22.6]{bib:Borel}, (i) $\tT = \pi^{-1}(T)$ is a maximal torus of~$\tG$, (ii) the maximal split subtorus
  $\tS \subset \tT$ satisfies $\pi(\tS) = S$, (iii) $X^*(S) \INTO X^*(\tS)$ induces a bijection $\alpha \mapsto
  \widetilde\alpha = \alpha \circ \pi$ on relative roots, (iv) $U_{\widetilde\alpha}$ maps isomorphically to $U_\alpha$ for any
  $\alpha \in \Phi$. Let $\Theta \subset \Phi$ be the set of roots of~$(S,L)$.
  Since $\tL = \langle \tT, U_{\widetilde\alpha}: \alpha \in \Theta\rangle$, $\tN \cong \prod_{\Phi^+-\Theta^+}
  U_{\widetilde\alpha}$ (in any fixed order) and similarly for $L$, $N$, the map~$\pi$ induces
  \begin{equation*}
    1 \to R \to \tL \to L \to 1, \quad \tN \congto N.
  \end{equation*}
  As $R$ is connected, $H^1(\Gal(\k/k), R(\k)) = 0$ by Lang's theorem so that
  \begin{align*}
    \tG(k) \onto G(k),\quad \tL(k) \onto L(k), \quad \tN(k) \congto N(k).
  \end{align*}
  Thus $V$ is an irreducible representation of~$\tG(k)$ on which $R(k)$ acts trivially.
  The result now follows from the previous case.
\end{proof}

The following technical lemma is crucial in controlling the support of the image of the Satake map.  Let~$\Phi\nd$ denote the
set of non-divisible roots in~$\Phi$.  Recall that for any root $\beta \in \Phi\nd$ there is a root subgroup~$U_\beta$
over~$F$ whose Lie algebra is the sum of weight spaces for the positive multiples of~$\beta$. It extends to a smooth
$\O$-subgroup scheme of~$G$ (see \S\ref{sec:buildings}).

\begin{lm}\label{lm:root_subgroup}
  Let~$\alpha$ be a simple root \(so $\alpha \in \Phi\nd^+$\).
  \begin{enumerate}
  \item The product map $\prod_{\beta \in \Phi\nd^+, \beta \ne \alpha} U_\beta \to U$ is an isomorphism of $\O$-schemes onto a closed
    subgroup scheme~$U'$. It is normal in~$U$ and independent of the order of the factors in the product. The product
    map induces an isomorphism of $\O$-group schemes $U_\alpha \ltimes U' \to U$.
  \item 
    Suppose that $A$ is an abelian group and that $\phi : U(F)/U(\O) \to A$ is a function with finite support. Then
    \begin{equation*}
      \sum_{U(F)/U(\O)} \phi(u) = \sum_{u_\alpha \in U_{\alpha}(F)/U_{\alpha}(\O)} \sum_{u' \in U'(F)/U'(\O)} \phi(u_\alpha u').
    \end{equation*}
  \item 
    Suppose that $\lambda \in X_*(S)$ and $\alpha \in \Phi\nd$ are such that $\langle \lambda,\alpha\rangle > 1$. Let $t =
    \lambda(\varpi)$.  Suppose that $A$ is an abelian group of exponent~$p$. Suppose that $\psi : U_\alpha(F)/tU_\alpha(\O)t^{-1} \to A$
    is a function with finite support such that $\psi$ is left invariant under $\ker(U_\alpha(\O) \to U_\alpha(k))$.  Then
    \begin{equation*}
      \sum_{u_\alpha \in U_{\alpha}(F)/t U_{\alpha}(\O)t^{-1}} \psi(u_\alpha) = 0.
    \end{equation*}
  \end{enumerate}
\end{lm}

\begin{proof}
  We will prove (i) and~(iii) at the end of Section~\ref{sec:buildings}. Part~(ii) follows immediately from~(i).

  Note however that when $G$ is split, the proof is easier.  In that case there are $\O$-group isomorphisms $x_\alpha : \G_a
  \congto U_\alpha$ such that for $t \in T$, $tx_\alpha (a) t^{-1} = x_\alpha(\alpha(t)a)$ and such that for all $\alpha$,
  $\beta$ with $\alpha \ne -\beta$, $[x_\alpha(a),x_\beta(b)] = \prod_{i,j > 0} x_{i\alpha + j\beta} (c_{i,j}a^ib^j)$ (in
  some order) with $c_{i,j} \in \O$. See \cite[II.1.2]{bib:Jan-reps}. Then (iii) is obvious since $U_\alpha$ is abelian and
  $tU_\alpha(\O)t^{-1}$ is a proper subgroup of $\ker(U_\alpha(\O) \to U_\alpha(k))$ of $p$-power index. Part~(i) follows as in
  the general case but instead of Bruhat--Tits one can appeal to \cite[II.1.7]{bib:Jan-reps}.
\end{proof}

\begin{proof}[Proof of Theorem~\ref{thm:satake}]
  We will use the refined Cartan decomposition (Lemma~\ref{lm:cartan})
  \begin{equation*}
    G(F) = \coprod_{\lambda \in X_*(S)_-} K\lambda(\varpi)K.
  \end{equation*}

  \emph{Step 0}. Let $f \in \HH_G(V)$.
  Since $K$ is compact open in~$G(F)$, $f$ is supported on a finite number of cosets in~$G(F)/K$. By the
  Iwasawa decomposition (Lemma~\ref{lm:iwasawa}), $f$ is supported on a finite number of cosets in $B(F)/B(\O)$. Thus $\SS f$
  is supported on a finite number of cosets in $T(F)/T(\O)$ and for each $t \in T(F)$, the sum $\sum_{u \in U(F)/U(\O)}
  f(tu)|_{V^{U(k)}}$ is zero outside a finite number of terms.  As $T$ normalises~$U$, it follows that the image of $\sum_{u
    \in U(F)/U(\O)} f(tu)|_{V^{U(k)}}$ is contained in~$V^{U(k)}$.  It is clear that $\SS$ is $\k$-linear.

  \emph{Step 1}. Show that the space of functions in~$\HH_G(V)$ supported on any single double coset is one-dimensional. 
  The argument is analogous to \cite[Lemma 7]{bib:BL2}, but it requires a technical input from Bruhat--Tits theory.
  Suppose that $f \in \HH_G(V)$ is supported on the double coset~$KtK$ with $t = \lambda(\varpi)$ for some $\lambda \in
  X_*(S)_-$. Let $P_\lambda = L_\lambda \ltimes U_\lambda$ denote the parabolic subgroup of~$G_{/k}$ defined by
  $\lambda \in X_*(S)$ \cite[13.4.2, 15.4.4]{bib:Springer_LAG}. Note that $L_\lambda = L_{-\lambda}$ and that
  $P_{-\lambda} = L_\lambda \ltimes U_{-\lambda}$ is the opposite parabolic.
  It follows immediately from the definitions that the possible values for~$f(t)$ are all $\phi \in \End_{\k} V$
  such that
  \begin{equation*}
    \text{$k_1\phi = \phi k_2$ whenever $k_1$, $k_2 \in K$ and $k_1 t = t k_2$.}
  \end{equation*}
  Note that $k_1 \in K \cap tKt^{-1}$, $k_2 \in K \cap t^{-1}Kt$ and $k_1 = t k_2 t^{-1}$.
  Prop.~\ref{prop:building} implies that equivalently $\phi$ has to factor through an $L_\lambda(k)$-equivariant map
  $V_{U_\lambda(k)} \to V^{U_{-\lambda}(k)}$ and Lemma~\ref{lm:U-inv} shows that the space of such~$\phi$ is one-dimensional (Schur's lemma).

  Again by Lemma~\ref{lm:U-inv} there is a function in~$\HH_G(V)$ that is supported on~$KtK$ and maps~$t$ to the endomorphism
  \begin{equation}\label{eq:3}
    V \onto V_{U_\lambda(k)} \tocong V^{U_{-\lambda}(k)} \INTO V.
  \end{equation}
  We denote it by~$T_\lambda$. Obviously it is a projection.

  \emph{Step 2}. Let us verify that $\SS$ is a homomorphism. This imitates the classical argument. Suppose that
  $f_i : G(F) \to \End_{\k} V$ ($i = 1$, $2$) are elements of~$\HH_G(V)$. Let $v \in V^{U(k)}$. Then
  $\SS(f_1 * f_2)(t)v$ equals
  {\allowdisplaybreaks
  \begin{align*}
    & = \sum_{u \in U(F)/U(\O)} \sum_{g \in G/K} f_1(tu g) f_2(g^{-1}) v \\
    &= \sum_{u \in U(F)/U(\O)} \sum_{b \in B(F)/B(\O)} f_1(tu b) f_2(b^{-1}) v \\
    &= \sum_{u \in U(F)/U(\O)} \sum_{\tau \in T(F)/T(\O)} \sum_{\nu \in U(F)/U(\O)} f_1(tu \tau\nu) f_2(\nu^{-1} \tau^{-1}) v \\
    &= \sum_{\tau \in T(F)/T(\O)} \sum_{\nu \in U(F)/U(\O)} \sum_{u \in U(F)/U(\O)} f_1(t\tau\nu) f_2(\nu^{-1} \tau^{-1} u) v \\
    &= \sum_{\tau \in T(F)/T(\O)} \sum_{\nu \in U(F)/U(\O)} \sum_{u \in U(F)/U(\O)} f_1(t\tau\nu) f_2(\tau^{-1} u) v \\
    &= \sum_{\tau \in T(F)/T(\O)} \sum_{\nu \in U(F)/U(\O)} f_1(t\tau\nu) (\SS f_2)(\tau^{-1}) v \\
    &= \sum_{\tau \in T(F)/T(\O)} (\SS f_1)(t\tau) (\SS f_2)(\tau^{-1}) v \\
    &= (\SS f_1 * \SS f_2)(t)v.
  \end{align*}}%
  Note that when we sum over quotients, the summand does not depend on the representative chosen provided we respect the
  stated order of summation.  The first and the last three equalities come from the definitions, the second from the Iwasawa
  decomposition $G(F) = B(F)K$ (Lemma~\ref{lm:iwasawa}), and the third follows from the fact that $B = T \ltimes U$.
  For the fourth equality, we replaced $(\tau^{-1} u \tau) \nu$ by~$\nu$, and for the fifth, we replaced
  $(\tau \nu^{-1} \tau^{-1})u$ by~$u$.

  \emph{Step 3}. Show that $(\SS T_\lambda)(\mu(\varpi)) = 0$ for $\mu \in X_*(S)$ unless $\mu \ge_\R \lambda$ and that
  $(\SS T_\lambda)(\lambda(\varpi)) = 1$. The argument is the classical one. By Lemma~\ref{lm:vanishing},
  $K\lambda(\varpi)K \cap \mu(\varpi)U \ne \varnothing$ implies $\mu \ge_\R \lambda$ and
  $K\lambda(\varpi)K \cap \lambda(\varpi)U = \lambda(\varpi)U(\O)$.
  Since $U_{-\lambda}(k) \subset U(k)$ and $T_\lambda(\lambda(\varpi))$ is a projection onto $V^{U_{-\lambda}(k)}$, we see
  that $(\SS T_\lambda)(\lambda(\varpi)) = 1$.

  \emph{Step 4}. Show that $(\SS f)(\mu(\varpi)) = 0$ if $\langle \mu, \alpha\rangle > 1$ for some simple root~$\alpha$.
  Let $t' = \mu(\varpi)$. By Lemma~\ref{lm:root_subgroup}(i), (ii),
  $U = U_\alpha \ltimes U'$ for some normal $\O$-subgroup scheme~$U'$ and for $v \in V^{U(k)}$,
  \begin{align*}
    (\SS f)(t')v &= \sum_{u_\alpha \in U_\alpha(F)/U_\alpha(\O)} \sum_{u' \in U'(F)/U'(\O)} f(t'u_\alpha u')v \\
    &= \sum_{u_\alpha \in U_\alpha(F)/t'U_\alpha(\O)t'^{-1}} \left(\sum_{u' \in U'(F)/U'(\O)} f(u_\alpha t'u')v\right).
  \end{align*}
  By Lemma~\ref{lm:root_subgroup}(iii) this sum is zero since $\langle \mu, \alpha\rangle > 1$, since $\k$ is of characteristic~$p$, and since the function
  of~$u_\alpha$ defined by the expression in parentheses is left invariant under 
  $\ker(U_\alpha(\O) \to U_\alpha(k))$.

  \emph{Step 5}. Show that $(\SS T_\lambda)(\mu(\varpi)) = 0$ if $\mu \not\in X_*(S)_-$. Suppose this is not the case.  Let
  $M_\lambda = \{ \mu \in X_*(S) : (\SS T_\lambda)(\mu(\varpi)) \ne 0 \}$.  Note that this is a finite set by Step~0.
  Label the simple roots $(\alpha_i)_{i=1}^r$ so that
  $\langle \mu, \alpha_1\rangle > 0$ for some $\mu \in M_\lambda$. Define a homomorphism of abelian groups
  \begin{align*}
    o : X_*(S) &\to \Z^r \\
    \mu &\mapsto (\langle \mu, \alpha_i\rangle)_{i=1}^r.
  \end{align*}
  Note that this is injective on $M_\lambda$: if $o(\mu_1) = o(\mu_2)$ for $\mu_i \in M_\lambda$, then $\mu_1-\mu_2 \in X^*_F(G)^\perp$
  (as $\mu_i \ge_\R \lambda$ by Step~3) and $\mu_1-\mu_2 \in (\R\langle \Phi\rangle)^\perp$, so $\mu_1 = \mu_2$.
  Let~$\mu$ be the
  element of~$M_\lambda$ such that $o(\mu)$ is greatest in the lexicographic order of~$\Z^r$. In particular, $\langle \mu,
  \alpha_1\rangle > 0$. We show that $\SS(T_\lambda^2) = (\SS T_\lambda)^2$ is non-zero on $2\mu(\varpi)$. Consider
  \begin{equation*}
    (\SS T_\lambda)^2(2\mu(\varpi)) = \sum_{\mu' \in X_*(S)} \SS T_\lambda(\mu'(\varpi)) \SS T_\lambda((2\mu-\mu')(\varpi)).
  \end{equation*}
  If the term indexed by~$\mu'$ is non-zero then $\mu'$, $2\mu-\mu' \in M_\lambda$ and hence $o(\mu') \le o(\mu)$ and
  $o(2\mu-\mu') \le o(\mu)$.  But since the sum of these inequalities yields an equality, it follows easily that $\mu' =
  \mu$. So $(\SS T_\lambda)^2(2\mu(\varpi)) = ((\SS T_\lambda)(\mu(\varpi)))^2 \ne 0$.
  Since $\langle 2\mu, \alpha_1\rangle > 1$, we get a contradiction by Step~4 with $f = T_\lambda^2$.

  \emph{Step 6}. It remains to show that $\SS$ is injective and maps onto $\HH_T^-(V^{U(k)})$. This is again classical.  By
  Step~1, the $T_\lambda$ ($\lambda \in X_*(S)_-$) form a $\k$-basis of~$\HH_G(V)$ and by Lemma~\ref{lm:torus-cosets}, the
  $\tau_\mu$ ($\mu \in X_*(S)_-$) form a $\k$-basis of~$\HH_T^-(V^{U(k)})$.  By Step~3, we may write $\SS T_\lambda =
  \sum_{\mu \ge_{\R} \lambda} a_\lambda(\mu) \tau_\mu$ with $a_\lambda(\mu) \in \k$ and $a_\lambda(\lambda) = 1$.  Since
  $\{\mu \in X_*(S)_- : \mu \ge_\R \lambda\}$ is finite by Lemma~\ref{lm:partial-order}, the claims follow.
\end{proof}

Suppose now that $G$ is split and that $G'$ is simply connected. We give a sketch of a simpler proof that $\HH_G(V)$ is
commutative.  By \cite[II.1.16]{bib:Jan-reps} there is a ``transpose'' involution ${}^\tau : G \to G$ that induces the
identity on $T$. (When $G = \GL_n$, one can take the usual transpose map.) Let ${}^\tau V$ be the dual $\Hom_\k (V,\k)$ with
$G(k)$-action $(g\psi)(v) := \psi({}^\tau g \cdot v)$.  Since $G'$ is simply connected, $V$ extends to a representation of
the algebraic group~$G_{/\k}$. By using a weight space decomposition of $V$, it follows that $V$ and ${}^\tau V$ are
isomorphic as $G(k)$-representations \cite[II.2.12(2)]{bib:Jan-reps}. Fix a $G(k)$-linear isomorphism $\kappa : V \congto {}^\tau
V$.

An element $\varphi \in \End_\k V$ induces an endomorphism of ${}^\tau V$ and hence an endomorphism ${}^\tau \varphi \in \End_\k V$ by using $\kappa$.
Given $f \in \HH_G(V)$, we define $f^* : G \to \End_\k V$ by $f^*(g) := {}^\tau f({}^\tau g)$. It is easy to check that 
$f^* \in \HH_G(V)$ and that $f_1^* \ast f_2^* = (f_2 \ast f_1)^*$. It remains to show that ${}^*$ acts trivially or equivalently that
$T_\lambda^* = T_\lambda$ for all $\lambda \in X_*(S)_-$. As ${}^\tau$ preserves $K = G(\O)$ and $\lambda(\varpi)$, it follows
that $T_\lambda^*$ has the same support as $T_\lambda$. Moreover it is clear that $T_\lambda^*(\lambda(\varpi))$ is a linear projection.
Hence $T_\lambda^* = T_\lambda$.

\subsection{Comparison with the $p$-adic Satake map}

We will explain a compatibility result with the $p$-adic Satake isomorphism of Schneider--Teitelbaum \cite[\S 3]{bib:ST}.
It will be convenient to state it in a slightly different form.
To keep the notation simple, let us assume in this subsection that $G_{/F}$ is split (just as in~\cite{bib:ST}).

Let $E$ be the (absolutely) irreducible representation of~$G_{/F}$ of highest weight $\nu \in X^*(T)$.  Then $E^{U(\O)}$ is
the highest weight space of~$E$; in particular it is one-dimensional and $T(\O)$ acts on it via~$\nu$. (This is because 
$E^{U(\O)} \subset E^{\mathfrak u}$, where $\mathfrak u = \Lie U(\O) = \Lie U(F) = \Lie (U_{/F})$. But $E^{\mathfrak u} = E^U$ since
$U_{/F}$ is connected.) 
Consider the following $p$-adic Hecke algebra,
\begin{equation*}
  \tHH_G(E) = \End_{G(F)} (\ind_K^{G(F)} E),
\end{equation*}
which we again think of as algebra (under convolution) of functions $f : G(F) \to \End_F(E)$ with compact support such that
$f(k_1gk_2) = k_1f(g)k_2$ for all $k_1$, $k_2 \in K$ and $g \in G(F)$.

\begin{lm}[{\cite[Lemma 1.4]{bib:ST}}]\label{lm:twist-hecke-alg}
  The map
  \begin{equation*}
    \iota : \tHH_G(1) \to \tHH_G(E)
  \end{equation*}
  with $(\iota \phi)(g) = \phi(g)g \in \End_F(E)$ is an algebra isomorphism.
\end{lm}

The point is that for $f \in \tHH_G(E)$ and $g \in G(F)$, we have $g^{-1}f(g) \in \End_F(E)^{K \cap g^{-1}Kg} = \End_F(E)^{G} = F$,
by considering the action of the Lie algebra as above.
Note that the lemma crucially depends on~$E$ being a representation not just of~$K$ but of~$G(F)$, thus the analogue does not
work for the characteristic~$p$ Hecke algebras.

Fix a $K$-stable norm $||.||_E$ on~$E$ such that $||E||_E = |F|$. Equivalently this corresponds to a choice of
$K$-stable $\O$-lattice $E_0 \subset E$ given by $E_0 = \{x \in E: ||x||_E \le 1\}$. Then $\tHH_G(E)$
carries a sub-multiplicative sup-norm, where $\End_F(E)$ is given the operator norm with respect to~$||.||_E$.
Similarly we have the Hecke algebra $\tHH_T(E^{U(\O)})$, likewise equipped with a sup-norm.
The $p$-adic Satake map is then the following isometric isomorphism of normed $F$-algebras:
\begin{align*}
  \tSS : \tHH_G(E) &\congto \tHH_T(E^{U(\O)})^{W,*} \\
  f &\mapsto \left(t \mapsto \sum_{U(F)/U(\O)} f(tu)\Big|_{E^{U(\O)}}\right).
\end{align*}
To define the right-hand side, let $\delta : B(F) \to q^\Z \subset \R\s$ be the modulus character of the Borel.
(Note that our~$\delta$ is inverse to the one in \cite{bib:ST}.)
Then $\tHH_T(E^{U(\O)})^{W,*}$ is the subalgebra of those $\varphi \in \tHH_T(E^{U(\O)})$ such that
$\varphi \nu^{-1} \delta^{1/2} : T(F)/T(\O) \to \overline F$ is $W$-invariant. This condition does not depend on the
choice of square root of~$\delta$
(see ~\cite{bib:ST}, Example 2 in \S 2). To prove that $\tSS$ is an algebra isomorphism one reduces to the case $E = 1$ by
applying Lemma~\ref{lm:twist-hecke-alg} to both sides, in which case it is equivalent to the classical Satake isomorphism.
That $\tSS$ is an isometry follows from Lemma~\ref{lm:vanishing}. For details, see \cite[\S 3]{bib:ST}. Note that their map
$S_\nu : \tHH_G(1) \to F[X_*(S)]$ \cite[p.~653]{bib:ST} is related to the above one by $S_\nu(\psi) = \varpi^{\ord \nu}
\nu^{-1} \tSS(\iota\psi)$.

From now on suppose that $E_0$ is a $G_{/\O}$-stable $\O$-lattice and that $E_0 \otimes_\O \k$ contains $F(\nu)$, the irreducible
representation of~$G_{/\k}$ of highest weight~$\nu$, as \emph{subobject}. For example, we could take the dual Weyl module $E_0 = H^0_\O(\lambda)$,
in the notation of \cite[II.8.6(1)]{bib:Jan-reps} (see also \cite[II.8.8(1), II.2.4]{bib:Jan-reps}).
Suppose moreover that $\nu$ is $q$-restricted, i.e., that $0 \le \langle \nu, \alpha\dual\rangle < q$ for all simple
roots~$\alpha$. Then $F(\nu)$ is irreducible as representation of~$G(k)$ and we denote it by~$V$. (See the proof of
Lemma~\ref{lm:U-inv}.  If $(G_{/\k})'$ is not simply connected, this follows by a $z$-extension argument.)
Let $\tHH_G(E)_0 \subset \tHH_G(E)$
denote the elements of sup-norm at most 1.  In particular, $\im(f) \subset \End_\O(E_0)$ for $f \in \tHH_G(E)_0$ and we can
consider the reduction $\overline f : G(F) \to \End_\k(E_0 \otimes \k)$.  Similarly we have $\tHH_T(E^{U(\O)})_0^{W,*}
\subset \tHH_T(E^{U(\O)})^{W,*}$. By considering the weight space decomposition of~$E_0$, it is clear that $E_0 \cap
E^{U(\O)}$ reduces to $V^{U(k)} \subset E_0 \otimes_\O \k$.

\begin{prop}\label{prop:st-comparison}
  With the above notation, we have the following commutative diagram.
  \ifkuvio
  \[ \cellwidth=45pt \Diag
  \tHH_G(E)_0 & \rTo_\sim^{\tSS} & \tHH_T(E^{U(\O)})_0^{W,*} \\
  \dTo_\alpha  && \dTo^\beta  \dy{-3mm} \\
  \HH_G(V) & \rTo_\sim^{\SS} & \HH_T^-(V^{U(k)}) \\ 
  \endDiag \]
  \fi
  Here $(\alpha f)(g) = \overline{f(g)}|_V$ and $(\beta \varphi)(t) = \overline{\varphi(t)}$.
  The vertical maps are well defined and induce isomorphisms after base extending from~$\O$ to~$\k$.
\end{prop}

\begin{proof}
  For $\lambda \in X_*(S)_-$ consider $\tT_\lambda \in \tHH_G(E)$ defined by (i) $\supp \tT_\lambda = K\lambda(\varpi)K$ and
  (ii) $\tT_\lambda(\lambda(\varpi)) = \varpi^{-\langle \lambda,\nu\rangle} \lambda(\varpi)$. We claim that the $\tT_\lambda$
  form an $\O$-basis of~$\tHH_G(E)_0$ and that $\alpha(\tT_\lambda) = T_\lambda$. On the $\nu'$-weight space of~$E$, for
  $\nu' \le \nu$, $\varpi^{-\langle \lambda,\nu\rangle} \lambda(\varpi)$ acts as the scalar $\varpi^{\langle \lambda,\nu' -
    \nu\rangle}$. Thus $\overline{\tT_\lambda(\lambda(\varpi))}$ is the linear projection onto the $\nu'$-weight spaces of
  $E_0\otimes_k \k$ for the weights~$\nu'$ satisfying $\langle \lambda,\nu' - \nu\rangle = 0$. Thus it preserves any
  $G_{/\k}$-subrepresentation, in particular, $V$. By~\eqref{eq:3} and by the description of~$V^{U_{-\lambda}(k)}$ given in 
  the proof of Lemma~\ref{lm:U-inv}, the claim follows and we see that $\alpha$ is well defined.
  
  Similarly, for $\lambda \in X_*(S)_-$ consider $\ttau_\lambda \in \tHH_T(E^{U(\O)})^{W,*}$ defined by (i) $T^- \cap \supp
  \ttau_\lambda = \lambda(\varpi)T(\O)$ and (ii) $\ttau_\lambda(\lambda(\varpi)) = 1$. We claim that the~$\ttau_\lambda$
  form an $\O$-basis of $\tHH_T(E^{U(\O)})^{W,*}_0$ and that $\beta(\ttau_\lambda) = \tau_\lambda$.  Recall
  that $\delta^{1/2}(\mu(\varpi)) = q^{-\langle \mu,\rho\rangle}$ for $\mu \in X_*(S)$, where $\rho = \frac 12 \sum_{\Phi^+}
  \alpha$ \cite[(3.3)]{bib:Gross_Satake}. Thus for $\varphi \in \tHH_T(E^{U(\O)})^{W,*}$,
  \begin{equation*}
    \varphi(w(\lambda(\varpi))) = \varphi(\lambda(\varpi)) \varpi^{\langle w\lambda - \lambda,\nu\rangle}
    q^{\langle w\lambda - \lambda,\rho\rangle} \quad \forall w \in W.
  \end{equation*}
  Since $w\lambda \ge_\R \lambda$ and since the second exponent is positive if $w\lambda \ne \lambda$, it follows that
  $\supp(\overline\varphi) \subset T^-$ whenever $||\varphi|| \le 1$. By the same reasoning, $||\ttau_\lambda|| \le 1$.
  The claim follows and we see that $\beta$ is well defined.

  This completes the proof since the diagram obviously commutes.
\end{proof}

\begin{rk}
  Note that this argument yields another proof that $\im(\SS) \subset \HH_T^-(V^{U(k)})$ in case
  $V$ arises from a representation of~$G_{/\k}$ \(which does not always happen if $(G_{/\k})'$ is not simply connected\),
  after the surjectivity of the map~$\alpha$ has been established.
\end{rk}

\subsection{The $W$-regular case}

For the proof of Prop.~\ref{prop:satake_regular} we will need a lemma. Let~$\oPhi$ denote the set of absolute roots of~$G_{/\k}$ with
respect to $T_{/\k}$. Since $G_{/k}$ is quasi-split, $W$ is a subgroup of the absolute Weyl group~$\oW$ and the restriction homomorphism
$X^*(T_{/\k}) \onto X^*(S_{/\k})$ is $W$-equivariant. Moreover, $\oPhi$ maps onto~$\Phi$ under this map; in particular, $\Phi^+$
determines a system of positive roots $\oPhi^+$ in~$\oPhi$.

\newcommand{\ob}{\overline\beta}

\begin{lm}\label{lm:roots}\
  \begin{enumerate}
  \item Suppose that $\eta \in X^*(T_{/\k})_+$ and $w \in W$. There are simple reflections $s_i \in W$ such that
    \begin{equation}\label{eq:7}
      \eta \gneq s_1 \eta \gneq \dots \gneq s_l\cdots s_2s_1 \eta = w\eta.
    \end{equation}
  \item Suppose that $\eta \in X^*(T_{/\k})_+$ and that $\alpha \in \Phi$ is simple. If $\eta - s_\alpha \eta \ge 0$ then
    $\eta - s_\alpha$ is the sum of simple roots $\ob_i \in \oPhi$ such that $\ob_i|_S = \alpha$.
  \end{enumerate}
\end{lm}

\begin{proof}
(i) Let us write $w = s_l\cdots s_1$ as a reduced product of simple reflections in~$W$. We will show that
\begin{equation}\label{eq:5}
  \eta \ge s_1 \eta \ge \dots \ge s_l\cdots s_2s_1 \eta = w\eta,
\end{equation}
which implies~\eqref{eq:7} since every time there is an equality, the corresponding simple reflection~$s_i$ can be omitted.
We claim that $\ell_\oW(w) = \sum \ell_\oW(s_i)$, where $\ell_\oW$ denotes the length in~$\oW$. Once we know this, we are done: by
writing each~$s_i$ as a reduced product of simple reflections in~$\oW$ we are reduced to proving the analogue of~\eqref{eq:5}
in~$\oW$, where it is easy and well known.

Recall that the length of~$w$ in~$W$ (resp.\ $\oW$) equals the number of non-divisible positive roots~$\alpha$ in~$\Phi$
(resp.\ $\oPhi$) such that $w(\alpha) < 0$ (see, for example, \cite{bib:Bourbaki-Lie}, VI.1.6, Cor.~2). In particular, a
simple reflection $s_\alpha \in W$ stabilises $\Phi^+-\{\alpha\}$. Say $\alpha_i \in \Phi$ is the simple root
corresponding to $s_i \in W$. Since $w = s_l\cdots s_1$ is of length~$l$ in~$W$, it sends precisely the following~$l$
non-divisible positive roots of~$\Phi$ to a negative root: $\alpha_1$, $s_1\alpha_2$, \dots, $s_1\cdots s_{l-1}\alpha_l$.
Letting $A_i = \{\ob \in \oPhi^+ : \ob|_S \in \Z_{>0} \alpha_i \}$, we see that $w$ sends precisely the following positive
roots of~$\oPhi$ to a negative root: $A_1 \cup s_1A_2 \cup \cdots \cup s_1\cdots s_{l-1}A_l$. Clearly $\ell_\oW(s_i) = |A_i|$, which
implies the claim.

(ii) Write $\eta - s_\alpha \eta = \ob_1+\cdots+\ob_r$ with $\ob_i \in \oPhi$ simple. Now restrict to~$S$. On the left-hand
side we get an integer multiple of~$\alpha$ and on the right-hand side a sum of simple roots $\ob_i|_S$ in~$\Phi$.
Thus $\ob_i|_S = \alpha$ for all~$i$.
\end{proof}

\begin{proof}[Proof of Proposition~\ref{prop:satake_regular}]
  By Step~3 of the proof of Thm.~\ref{thm:satake}, we know that $(\SS T_\lambda)(\lambda(\varpi)) = 1$.
  It thus suffices to show that for any given $\mu \in X_*(S) - \{\lambda\}$,
  each term in the sum defining $(\SS T_\lambda)(\mu(\varpi))$ vanishes. 

  Let $t' = \mu(\varpi)$, $t = \lambda(\varpi)$.
  Choose $0 \ne v \in V^{U(k)}$.

  \emph{Step 1.} We will show that if $T_\lambda(g) v \ne 0$ then $g \in Kt I$, where $I = \red^{-1} (B(k))$ is an Iwahori
  subgroup. Let $W_\lambda \le W$ be the Weyl group of~$(S_{/k},L_\lambda)$ (generated by simple reflections associated to simple
  roots $\alpha \in \Phi$ with $\langle \lambda,\alpha \rangle = 0$). For each $w \in W$ choose a representative $\dot w
  \in N(S)(k)$ and a lift of it, $\dot w \in G(\O) = K$. Then
  \begin{equation*}
    G(k) = \coprod_{W_\lambda \backslash W} P_\lambda(k) \dot w B(k)
  \end{equation*}
  by~\cite[21.16(3)]{bib:Borel}. By Prop.~\ref{prop:building},
  \begin{equation*}
    K = \coprod_{W_\lambda \backslash W} (K \cap t^{-1} K t) \dot w I
  \end{equation*}
  and thus
  \begin{equation*}
    KtK = \coprod_{W_\lambda \backslash W} Kt \dot w I.
  \end{equation*}
  So if $T_\lambda(g) v \ne 0$ then $g = k t \dot w i$ for some $k \in K$, $w \in W$, $i \in I$. Thus $T_\lambda(t)\dot w v
  \ne 0$. We will show that $w \in W_\lambda$.  Recalling the definition of~$T_\lambda$ \eqref{eq:3}, we may by the proof of
  Lemma~\ref{lm:U-inv} reduce to the case that $(G_{/\k})'$ is simply connected. (The lifted Levi equals
  $L_{\widetilde\lambda}$ for any lift $\widetilde\lambda \in X_*(\tS)$ of~$\lambda$.  We can lift $\dot w$ since the Weyl
  groups of~$(S_{/k},G_{/k})$ and its lift $(\tS,\tG)$ are naturally identified by \cite[22.6]{bib:Borel}.)  Since
  $(G_{/\k})'$ is simply connected, there is a $q$-restricted weight $\nu \in X^*(T)_+$ such that $V \cong F(\nu)$ as
  $G(k)$-representations. But we saw in the proof of Lemma~\ref{lm:U-inv} that $T_\lambda(t)$ is the projection onto the
  weight spaces for $\nu' \in X^*(T)$ such that $\nu-\nu'$ is a sum of simple roots of~$(T_{/\k},L_{\lambda/\k})$, i.e., a
  sum of simple roots $\ob \in \oPhi$ such that $\langle \ob,\lambda\rangle = 0$.  Since $T_\lambda(t)\dot w v \ne 0$, it
  follows that $\nu -w\nu$ is a sum of simple roots $\ob \in \oPhi$ such that $\langle \ob,\lambda\rangle = 0$.

  By Lemma~\ref{lm:roots}(i) there are simple reflections $s_i \in W$ corresponding to simple roots $\alpha_i \in \Phi$ such
  that
  \begin{equation*}
    \nu \gneq s_1\nu \gneq s_{2}s_1\nu \gneq \cdots \gneq s_l\cdots s_1 \nu = w\nu.
  \end{equation*}
  By Lemma~\ref{lm:roots}(ii) the $i$-th and $(i+1)$-st term in this sequence differ by a sum of simple roots $\ob_{ij} \in \oPhi$
  such that $\ob_{ij}|_S = \alpha_i$. Thus $\langle\alpha_i,\lambda\rangle = \langle \ob_{ij},\lambda\rangle = 0$.
  It follows that $s_i \in W_\lambda$ for all~$i$. Since $V$ is $W$-regular we see that $w =
  s_l\cdots s_1\in W_\lambda$ and $g \in K t I$.

  (We remark that we only actually used that $\Stab_W(\nu) \subset W_\lambda$.)
  
  \emph{Step 2.} We show that $K t I \cap t'U(F) = \varnothing$. Suppose not. We use the Iwahori decomposition
  \begin{equation*}
    I = (I \cap \oU(F))(I \cap T(F))(I \cap U(F)),
  \end{equation*}
  where $\oU$ is the unipotent radical of the opposite Borel (Lemma~\ref{lm:iwahori}). Since $t$ contracts $I \cap \oU(F)$ we find
  that $t I t^{-1} \subset I U(F)$.  Thus
  \begin{equation*}
    \varnothing \ne (K t I \cap t'U(F))t^{-1} \subset K U(F) \cap t't^{-1}U(F).
  \end{equation*}
  Therefore $K \cap t't^{-1} U(F) \ne \varnothing$ and so $t't^{-1} \in T(\O)$, which contradicts that $\mu \ne \lambda$.
\end{proof}

\subsection{Satake parameters}

\begin{proof}[Proof of Corollary~\ref{cor:satake_params}]
  By Cor.~\ref{cor:hecke_struct} we need to classify algebra homomorphisms $\theta : \k[X_*(S)_-]
  \to \k$, i.e., monoid homomorphisms $X_*(S)_- \to \k$, where $\k$ is considered with its multiplicative structure.
  Then $M := \theta^{-1} (\k\s)$ satisfies
  \begin{equation}\label{eq:9}
    \lambda_1 + \lambda_2 \in M \iff \lambda_1 \in M \text{\ and\ } \lambda_2 \in M.
  \end{equation}
  Let $X_*(S)_0 := \{\lambda \in X_*(S) : \langle \lambda,\alpha \rangle = 0 \ \forall \alpha \in \Phi\}$. Since this is a subgroup
  of~$X_*(S)_-$, $X_*(S)_0 \subset M$.
  For $\delta \in \Delta$ choose $\lambda_\delta \in X_*(S)_-$ such that
  $\langle \lambda_\delta, \delta'\rangle$ is zero if $\delta' \in \Delta-\{\delta\}$ and negative if $\delta' = \delta$.
  
  We claim that $M = J^\perp \cap X_*(S)_-$ (a ``facet'' of~$X_*(S)_-$), where $J = \{ \delta : \lambda_\delta \not\in M \}$.
  (Note that $J$ is independent of the choice of the $\lambda_\delta$, since $X_*(S)_0 \subset M$.)
  Suppose that $\lambda \in X_*(S)_-$.  Then there is an $n \in \Z_{>0}$ such that $n\lambda = \sum n_\delta \lambda_\delta +
  \lambda_0$ for some $n_\delta \in \Z_{\ge 0}$ and some $\lambda_0 \in X_*(S)_0$. Then from~(\ref{eq:9}) we see that
  $\lambda \in M$ iff $n_\delta \ne 0$ implies $\delta \not\in J$ iff $\lambda \in J^\perp$.
  
  Next we show that the subgroup of~$X_*(S)$ generated by~$M$ equals $J^\perp$. One inclusion being obvious,
  suppose that $\lambda \in J^\perp$. Then $\lambda + n\sum_{\delta \not\in J} \lambda_\delta \in X_*(S)_-$
  (and hence in $M =J^\perp \cap X_*(S)_-$) for some $n \in \Z_{>0}$, which implies that $\lambda$ is in the subgroup generated by~$M$.

  As $\k\s$ is a group, $\theta|_M$ extends uniquely to a group homomorphism $\tilde\theta : J^\perp \to \k\s$. Taking
  character groups in the exact sequence defining $\hS_J$, we find that $X^*(\hS_J) = J^\perp$. Thus $\tilde\theta$ corresponds
  to an element of $X_*(\hS_J) \tens \k\s \cong \hS_J(\k)$.

  All pairs $(J,s_J)$ with $s_J \in \hS_J(\k)$ are obtained in this way, because $J^\perp \cap X_*(S)_-$ satisfies~(\ref{eq:9}),
  which allows us to extend a homomorphism $J^\perp \cap X_*(S)_- \to \k$ by zero to a monoid homomorphism $X_*(S)_- \to \k$.
\end{proof}

\section{Buildings arguments}\label{sec:buildings}

\renewcommand{\B}{\mathfrak B}
\newcommand{\hB}{\widehat \B}
\renewcommand{\G}{\mathfrak G}
\newcommand{\Gx}{\G_x^0}
\newcommand{\GF}{\G_F^0}
\newcommand{\oG}{\overline{\G}}
\newcommand{\oGx}{\oG_x^0}
\newcommand{\hGx}{\widehat \G_x}
\newcommand{\hGF}{\widehat \G_F}
\newcommand{\hGO}{\widehat \G_\Omega}
\newcommand{\II}{\mathcal I}
\newcommand{\tII}{\widetilde\II}
\newcommand{\K}{\underline K}
\newcommand{\oK}{\overline K}
\newcommand{\NO}{N^1_\Omega}
\newcommand{\NF}{N^1_F}
\newcommand{\hN}{\widehat N^1}
\newcommand{\hNO}{\widehat N^1_\Omega}
\newcommand{\hNF}{\widehat N^1_F}
\newcommand{\tR}{\widetilde \R}
\newcommand{\oS}{\overline{\S}}
\renewcommand{\T}{\mathfrak T}
\newcommand{\oT}{\overline{\T}}
\newcommand{\U}{\mathfrak U}
\newcommand{\UO}{U_\Omega}
\renewcommand{\oU}{\overline{\U}}
\newcommand{\V}{\mathfrak V}
\newcommand{\hV}{\widehat V}
\newcommand{\vW}{{}^v W}
\newcommand{\hW}{\widehat W}

\newcommand{\vnu}{{}^v \nu}

The main goal of this section is to prove Prop.~\ref{prop:building} and Lemma~\ref{lm:root_subgroup}. We also justify
some basic results about unramified groups using the work of Bruhat--Tits \cite{bib:BT1}, \cite{bib:BT2}.
Although most of these are well known, we could not find a good reference for the proofs. 

References to chapters~I \cite{bib:BT1} and~II \cite{bib:BT2} of Bruhat--Tits will be given in the form I.4.4.4, II.5.1.40
(for example).

We will keep as much as possible with the notation of Bruhat--Tits. In particular $K$ now denotes the $p$-adic field
and~$\oK$ its residue field, $N$ denotes~$N(S)$, $Z$ denotes the centraliser $Z(S)$ of $S$ in $G$, and~$\vW$ the Weyl group.
Group schemes over~$\O$ are denoted by fraktur letters ($\G$, $\T$, \dots), their generic fibres by the corresponding roman
letters ($G$, $T$,\dots) and their special fibres are overlined ($\oG$, $\oT$,\dots).  Note that ``fixer'' is a synonym for
``pointwise stabiliser.'' An $\O$-group scheme is said to be \emph{connected} if its two fibres are connected. The connected
component of a smooth $\O$-group scheme is defined fibrewise (II.1.2.12).  As in \S\ref{sec:proofs} we are assuming that the
valuation $\ord_K$ surjects onto the integers.

Let~$\II$ denote the \emph{reduced} building of~$G$. The general construction of I.6 and I.7 produces~$\II$ starting with
a valuation of the ``root datum'' $(T(K),(U_a(K))_a)$. Such a valuation is constructed for quasi-split groups by descent from the split case (II.4.2)
and in general by \'etale descent from the quasi-split case (II.5.1). The apartment~$A$ of~$S$ is an affine space under the vector space~$V$
which is the quotient of~$X_*(S)_\R$ dual to $\R\langle \Phi\rangle \subset X^*(S)_\R$. 

\renewcommand{\S}{\mathfrak S}

\begin{lm}
  Suppose that $\G$ is a smooth $\O$-group scheme with generic fibre~$G$. Then $\G \times \oK$ is reductive if and only if
  $\G \cong \Gx$ for some hyperspecial point~$x$. In this case $G$ is unramified and $\G \times \oK$ is connected.
\end{lm}

Recall that a point $x \in \II$ is \emph{hyperspecial} if $G$ splits over~$K\nr$ and $x$ is a special point inside the building
of $G \times K\nr$ \cite[1.10]{bib:Tits}.

\begin{proof}
  The first statement is II.5.1.40. (Note that in II.5, the superscript~$\natural$ refers to the objects over the base field; the other objects live
  over the strict henselisation of the base field.) 
  
  Let us show that $G$ is quasi-split. Without loss of generality, assume that $x$ lies in the apartment of~$S$.  The
  canonical extension~$\S$ of~$S$ (the split torus over~$\O$ with generic fibre~$S$) is a closed subscheme of~$\Gx$ and its
  reduction~$\oS$ is a maximal $\oK$-split torus in~$\oGx$ (II.5.1.11). The Lie algebra $\Lie \Gx$ is a free $\O$-module of
  finite rank (as $\Gx$ is a smooth group scheme) and we can consider its decomposition under~$\S$. Note the character groups
  $X^*(\S)$, $X^*(\S_{\oK})$, $X^*(S)$ are naturally isomorphic. Since $\oK$ is a finite field, $\oGx$ is quasi-split and
  \begin{equation*}
    \rank \oGx = \dim_{\oK} (\Lie \oGx)^{\oS = 1} = \dim_K (\Lie G)^{S = 1} = \dim Z \ge \rank G.
  \end{equation*}
  (Here ``$\rank$'' denotes the absolute rank of an algebraic group.) On the other hand, any split torus in the special fibre
  of $\Gx\times \O\nr$ can be lifted to a split torus in the generic fibre, as explained in the proof of II.4.6.4, so that
  $\rank \oGx \le \rank G$.  Thus equality holds, so $Z$ is a maximal torus of~$G$, i.e., $G$ is quasi-split.

  The connectedness of $\G \times \oK$ follows by base change to the strict henselisation and II.4.6.22.
\end{proof}

Assume from now on that $\G \cong \Gx$ for some hyperspecial point~$x$. Then $\K := \Gx(\O)$ is a hyperspecial maximal compact subgroup of~$G(K)$.

Let us summarise some results in II.4.6 on the structure of $\G \cong \Gx$.  Fix an apartment~$A$ of~$\II$ containing~$x$.
Let $S$ be the corresponding maximal split torus of~$G$ and let $T = Z$ (a maximal torus, since $G$ is quasi-split). Let~$\Phi$
be the set of roots of~$(G,S)$ and let $\Phi^\red$ denote the subset of non-divisible roots.
For $a \in \Phi$ let~$U_a$ denote the corresponding root subgroup.  In particular,
$U_{2a} \subset U_a$ whenever $\{a, 2a\} \subset \Phi$.  Fix a Steinberg--Chevalley valuation $\varphi = (\varphi_a)_{a \in
  \Phi}$ of the ``root datum'' $(T(K),(U_a(K))_{a\in \Phi})$, as constructed in II.4.1, II.4.2.  Here $\varphi_a : U_a(K) \to \R
\cup \{\infty\}$. It yields a filtration of each root subgroup: $U_{a,k} = \{ u \in U_a(K): \varphi_a(u) \ge k \}$
(II.4.3.1(1)).  Let $\Gamma_a = \varphi_a(U_a-\{1\})$ and $\Gamma'_a = \{\varphi_a(u) : u \in U_a -\{1\}, \varphi_a(u) = \max
\varphi_a(u U_{2a})\} \subset \Gamma_a$; these are discrete subsets of~$\R$.

By II.4.4.18 there are smooth prolongations $\S$ of~$S$ (the split torus over~$\O$ with generic fibre~$S$) and $\T$ of~$T$ (denoted
there by~$\T^R$). Then $\S$ is a closed subgroup scheme of~$\T$.

\begin{lm}\label{lm:conn_torus}
  $\T$ is connected \(i.e., its special fibre is connected\).
\end{lm}

\begin{proof}
Let~$K\nr$ be the maximal unramified extension of~$K$ with ring of integers $\O\nr$. Since $T \times K\nr$ is split, it has
a canonical prolongation $\T\nr$ to $\O\nr$ (the split torus over~$\O\nr$ with generic fibre $T \times K\nr$).
As remarked in II.5.1.9 (top of p.~149), $\T\nr$ descends to the torus~$\T$ defined in II.4.4. Since $\T\nr$ is connected, this completes the proof.
To justify that remark in II.5.1.9, one uses the last item in II.4.4.12(i) and the fact that $\T\nr$ is \'etoff\'e (II.1.7)
to see that $\O[\T] = \{ f \in K[T] : f(\T\nr(\O\nr)) \subset \O\nr \} = \O[\T']$, where $\T'$ is the torus descended from $\T\nr$.
\end{proof}

From II.4.6.4 it follows that $\oS$ is a maximal split torus of~$\oGx$ and that $\oT$ is the centraliser of~$\oS$ (a maximal
torus, as $\oGx \times \oK$ is quasi-split). By considering the Lie algebra of~$\Gx$, we see that the root systems of~$(S,G)$
and $(\oS, \oGx)$ are naturally identified.

Recall that $\Gx$ is the smooth $\O$-group scheme $\G_{f}^0$ with generic fibre~$G$ associated to the optimal, quasi-concave
function $f : \Phi \to \R$,
\[ f(a) = \min\{ k \in \Gamma'_a : a(x-\varphi) + k \ge 0 \}.  \]
(See II.4.6.26.) For all non-divisible roots $a \in \Phi$, there is a smooth $\O$-group scheme $\U_{f,a}$ with generic fibre~$U_a$ (II.4.5)
which we denote by~$\U_{x,a}$.  It is a closed subgroup scheme of~$\Gx$ and $\oU_{x,a}$ is the root subgroup of~$a$ in~$\oGx$
(II.4.6.4). The product map $\prod_a \U_{x,a} \to \Gx$, where $a$ runs over all positive, non-divisible roots in any order,
is an isomorphism onto a closed subgroup scheme~$\U^+$ (II.4.6.2). Let~$U^+$ denote its generic fibre.
$\T$ normalises each $\U_{x,a}$ (II.4.4.19) and the
product map yields an isomorphism of the semidirect product $\T \ltimes \U^+$ onto a closed subgroup scheme
of~$\Gx$ whose fibres are the Borel subgroups associated to $\Phi^+$. (Note that this is stated in II.3.8.2 only for a group
scheme whose connected component is~$\Gx$, but this implies the assertion here: the scheme $\T \times \U^+$ is connected as it is the
product of connected group schemes \cite[Exp.\ VI$_A$, Lemme 2.1.2]{bib:SGA3}.)

\begin{lm}\label{lm:fixer}
  Suppose that $F$ is a facet of~$A$ whose closure contains the hyperspecial point~$x$. Then $\hGF = \GF$.
  In particular, $\hGx = \Gx$.
\end{lm}

Note that $\hGO$ (II.4.6.26) equals $\mathcal G_{\mathrm{pr}_{\mathrm ss}^{-1}(\Omega)}$ in the notation of~\cite[3.4]{bib:Tits}.

\begin{proof}
  First we show $\hGF = \G_F$ by showing that $\hNF = \NF$ (II.4.6.26).
  Let $G(K)^1 = \{ g \in G(K) : \ord_K(\chi(g)) = 0 \ \forall \chi \in X_K^*(G) \}$.
  Note that $\ker \nu \cap G(K)^1 = H^1 \subset \NF \subset
  \hNF$ (II.4.6.3) so it suffices to show that $\nu(\NF) = \nu(\hNF)$. Identify $A$ and~$V$ using the special point~$x$ as origin. Then
  $\hNF$ is identified with a subgroup of~$\vW$, namely the subgroup of elements fixing~$F$. It is generated by those basic
  reflections~$r_a$ of~$\vW$ 
  such that $F$ is contained in the hyperplane through~$x$ which is defined by $a \in \Phi$. But
  I.7.1.3 shows that $\nu(\NF)$ has the same description. (The point is that $\Gamma_a = \Gamma'_a \cup \frac 12 \Gamma_{2a}$
  (I.6.2.1) and that $r_a = r_{2a}$.)

  Finally, $\G_F = \GF$ since $\T$ is connected (II.4.6.2). (This is the only part that uses that $x$ is hyperspecial, not
  just special.)
\end{proof}

\begin{lm}[Iwasawa decomposition {\cite[3.3.2]{bib:Tits}}]\label{lm:iwasawa}
  \[ G(K) = T(K)U^+(K) \K\]
\end{lm}

\begin{proof}
  We use the description of the building in terms of an affine Tits system. Associated to the valuation~$\varphi$ of the
  ``root datum'' $\big(T(K),(U_a(K))\big)$, we have the apartment~$A$, the set of affine roots $\alpha_{a,k}$ ($a \in \Phi$,
  $k \in \Gamma'_a$), the affine Weyl group~$W$ generated by the set of reflections in the boundary hyperplanes of the affine
  roots and $\nu : N(K) \to \Aff(A)$ giving the action on the apartment with kernel~$H$ (I.6.2). Let $N' = \nu^{-1}(W)$, $T'
  = N' \cap T(K)$, and $G' = \langle N', U_a(K)\rangle_{a \in \Phi}$. Fix a chamber $C \subset A$. Let $B = HU_C$ and let
  $\mathbf{S}$ be the set of reflections in the walls of $C$. By I.6.5, $(G',B,N',\mathbf S)$ is a saturated Tits system of
  affine type such that the inclusion $G' \to G(K)$ is $(B,N')$-adapted of connected type and such that the condition $G' =
  \B N' B$ in I.4.4(1) holds with $\B = T'U^+(K)$.

  Then $\II$ is naturally isomorphic with the building constructed out of this Tits system, whose facets are the
  ``parahoric'' subgroups of $(G',B,N',\mathbf S)$ (I.2, I.7.4.2). 
  Let~$\K'$ be the fixer of~$x$ in~$G(K)$, so that $\K = \K' \cap G(K)^1$ by Lemma~\ref{lm:fixer}.  By I.4.4.5, $\K' =
  (\nu^{-1}(\hV) \cap \K')\K$, where $\hV$ consists of the translations in $\hW = \nu(N(K))$.  As $x$ is special, $\K'$ is a
  good maximally bounded subgroup of~$G(K)$ (I.4.4.6(i)) so that $G(K) = \hB \K' = \hB (\nu^{-1}(\hV) \cap \K')\K$.  The
  result follows by using that $\hB = \nu^{-1}(\hV) \B$ (I.4.1.5) and $\nu^{-1}(\hV) = T(K)$ (I.6.2.10(i), I.6.1.11(ii)).
\end{proof}

\begin{lm}[Cartan decomposition {\cite[3.3.3]{bib:Tits}}]\label{lm:cartan}
  \[ G(K) = \coprod_{\lambda \in X_*(S)_-} \K \lambda(\varpi) \K \]
\end{lm}

\begin{proof}
  We keep the notation of the previous proof. Let~$D$ be the ``Weyl'' chamber in~$V$ corresponding to $\Phi^+$ and let
  $\hV_D = \hV \cap \overline D$.  By I.4.4.3(2), $G(K) = \K' \nu^{-1}(\hV_D) \K'$ and the set of double cosets biject with~$\hV_D$.
  Since $\nu^{-1}(\hV) \cap \K' = \ker \nu = H$, $\K' = H\K$ and $G(K) = \K \nu^{-1}(\hV_D) \K$.  Besides, $G(K)^1
  \lhd G(K)$ and $H \subset T(K)$. Using these facts it is easy to see that for $t_1$, $t_2 \in \nu^{-1}(\hV_D) \subset
  T(K)$, $\K t_1\K = \K t_2\K$ if and only if $t_1t_2^{-1} \in H \cap G(K)^1 = \ker \nu^1$ where $\nu^1$ is the action map of~$N(K)$
  on the extended apartment (II.4.2.16). It follows that the set of double cosets $\K\backslash G(K)/\K$ bijects with
  $\hV^1_D = \nu^1(\nu^{-1}(\hV_D))$ (the analogue of~$\hV_D$ for the extended building).

  By I.4.2.16(3), $\langle \nu^1(t),c\rangle = -(\ord_K \circ c)(t)$ for $t \in T(K)$ and $c \in X^*_K(T)_\R = X^*(S)_\R$. By Lemma~\ref{lm:torus-cosets},
  $\nu^1(t) = -\zeta(t)$ where $\zeta : T(K) \to X_*(S)$ was defined there. The result follows from that lemma.
\end{proof}

\begin{lm}\label{lm:vanishing}
  Suppose that $\lambda \in X^*(S)_-$, $\lambda' \in X^*(S)$.

  \begin{enumerate}
  \item If $\K \lambda(\varpi) \K \cap \lambda'(\varpi)U^+(K) \ne \varnothing$ then $\lambda' \ge_{\R} \lambda$.
  \item $\K \lambda(\varpi) \K \cap \lambda(\varpi)U^+(K) = \lambda(\varpi)\U^+(\O)$.
  \end{enumerate}
\end{lm}

Note that (i) is claimed without proof in \cite[p.~148]{bib:Cartier}.

\begin{proof}
  We keep the notation of Lemmas~\ref{lm:iwasawa} and~\ref{lm:cartan}.
  
  (i) For $t \in T(K)$ and $a \in \Phi$ we have $\langle \nu(t), a \rangle = -(\ord_K \circ a)(t)$ (I.4.2.7(3)).  Applying
  this with $t = (\lambda-\lambda')(\varpi)$, we see that the image of $\lambda'-\lambda \in X_*(S)_\R$ in the quotient space
  $V \cong (\R\langle \Phi\rangle)^*$ is $\nu((\lambda-\lambda')(\varpi))$.

  Suppose that $\K \lambda(\varpi) \K \cap
  \lambda'(\varpi)U^+(K) \ne \varnothing$.  On the one hand, $\K' \lambda(\varpi) \K' \cap HU^+(K) \lambda'(\varpi)\K' \ne
  \varnothing$. This implies that $\nu((\lambda-\lambda')(\varpi)) \ge_D 0$ by I.4.4.4(i), i.e., $\bar\lambda' \ge_\R \bar\lambda$. 
  (Note that $\hB^0 = H\B^0 = HU^+(K)$ by I.4.1.5 and I.6.5.)
  On the other hand, as $\K$ and $U^+(K)$ are contained in~$G(K)^1$,
  $(\lambda'-\lambda)(\varpi) \in G(K)^1$, i.e., $\lambda'-\lambda \in X^*_K(G)_\R^\perp$. The assertion follows from the
  definition of~$\le_\R$ (Def.~\ref{df:partial-order}).

  (ii) Note that the left-hand side is contained in
  \begin{multline*}
    \big(\K' \lambda(\varpi) \K' \cap HU^+(K)\lambda(\varpi)\K'\big) \cap \lambda(\varpi)U^+(K) 
    = \lambda(\varpi)\K' \cap \lambda(\varpi)U^+(K)
  \end{multline*}
  by I.4.4.4(ii). As $U^+(K) \subset G(K)^1$, this is contained in $\lambda(\varpi)\K \cap
  \lambda(\varpi)U^+(K) = \lambda(\varpi)\U^+(\O)$.  The opposite containment is obvious.
\end{proof}

\begin{lm}\label{lm:hyperspecial}
  If $y \in A$ is hyperspecial then $a(\varphi-y) \in \Gamma'_a$ for all $a \in \Phi$.
\end{lm}

\begin{proof}
  We consider the $G(K)$-equivariant injection of buildings $j : \II \to \tII$, where $\tII$ is the building of~$G$
  over~$K\nr$ (II.5.1.24), or even just the restriction of~$j$ to apartments $A \to \widetilde A$ corresponding to~$S$
  (resp.~$T$).  Let~$\widetilde\Phi$ denote the set of roots of~$(T,G)$. For $a \in \Phi$ let us say that an ``$a$-wall'' is
  the boundary of an affine root defined by~$a$ in~$A$.  Similarly we have the notion of an ``$\widetilde a$-wall'' in~$\widetilde A$
  for $\widetilde a \in \widetilde\Phi$.  By II.5.1.20, the affine roots in~$A$ are precisely the
  intersections with~$A$ of the affine roots in~$\widetilde A$.

As $y$ is hyperspecial, for each $\widetilde a \in \widetilde\Phi$ there is an $\widetilde a$-wall passing through~$j(y)$.
By intersecting with~$A$, we see that there is an $a$-wall passing through~$y$ for each $a \in \Phi$. Since the affine roots
in~$A$ are defined to be the $\alpha_{a,k} = \{ z \in A : a(z-\varphi) + k \ge 0\}$ for $a \in \Phi$ and $k \in \Gamma'_a$,
the lemma follows.
\end{proof}

As in the proof of Thm.~\ref{thm:satake} we denote by $P_\lambda = L_\lambda \ltimes U_\lambda$ the parabolic subgroup
of $\G \times \oK$ determined by $\lambda \in X_*(\oS) = X_*(S)$.

\begin{prop}\label{prop:building}
  Suppose that $\lambda \in X_*(S)$. Let $t = \lambda(\varpi) \in S(K)$ and let $\red : \G(\O) \to \G(\oK)$ denote the reduction map. Then
  \begin{equation*}
    \red(\G(\O) \cap t^{-1}\G(\O)t) = P_\lambda(\oK).
  \end{equation*}
  Moreover,
  \begin{multline}\label{eq:10}
  \{ (\red(g), \red(tgt^{-1})) : g \in \G(\O) \cap t^{-1}\G(\O)t \} \\
  = \{(g_+,g_-) \in P_\lambda(\oK) \times P_{-\lambda}(\oK) : [g_+] = [g_-] \in L_\lambda(\oK) = L_{-\lambda}(\oK)\},
  \end{multline}
  where $[\cdot]$ denotes the projection to the Levi subgroup.
\end{prop}

Note that this is actually obvious when $\G = {\GL_n}$.

\begin{proof}
  Let $\Omega = \{x, t^{-1} x\} \subset A$. By Lemma~\ref{lm:fixer}, $\G(\O) \cap t^{-1}\G(\O)t$ is the fixer of~$\Omega$ in~$G(K)^1$.
  Thus it equals $\hNO \UO$ by I.7.4.4, II.4.6.26.
  
  Let us first show that $\nu(\hNO)$ is naturally isomorphic to $\vW_\lambda = \{ w \in \vW : w\lambda = \lambda\}$. As $x$ is
  special, $\nu(\hN_x)$ is isomorphic to~$\vW$ via the forgetful map $\Aff(A) \to \GL(V)$ (I.6.2.10).
  Suppose $n \in \hN_x$ and let $w = \vnu(n)$. Then $n \in \hNO$ iff $w$ fixes $x - t^{-1}x \in V$. By II.4.2.7(3), $\langle
  x-t^{-1}x,a\rangle = -\ord_K(a(t))$ for $a \in \Phi$. So $w$ fixes $x - t^{-1}x$ iff $\lambda - w\lambda \in \langle \Phi
  \rangle^\perp$. But for all $w \in \vW$, $\lambda -w\lambda \in X^*_K(G)^\perp$.  Thus $\lambda - w\lambda
  \in \langle \Phi \rangle^\perp$ is equivalent to $\lambda = w\lambda$.
  
  Next we show that $\red(\hNO)$ equals the $\oK$-points of the normaliser of~$\oS$ in~$L_\lambda$. Note that $\hNO \supset
  G(K)^1 \cap \ker \nu = \T(\O)$ (II.4.6.3(3)).  Also, $\hNO \subset N(K) \cap \G(\O)$. If $n \in N(K) \cap \G(\O)$ then, by considering
  the $\O\nr$-points of~$\S$, we see that $\red(n) \in N(\oS)$ and that $\red(n)$ induces the same Weyl element on $X_*(\oS)$ as~$n$
  on~$X_*(S)$. From the previous paragraph, $\red(\hNO)/\oT(\oK) \cong \vW_\lambda$, which is precisely the Weyl
  group of~$\oS$ in~$L_\lambda$.

  To determine $\red(\UO)$, let us compute $f_\Omega : \Phi \to \R$. By II.4.6.26, for $a \in \Phi$,
  \begin{align*}
    f_\Omega(a) & = f_x(a) + \max(a(x-t^{-1}x),0) \\
    &= f_x(a) + \max(-\langle a,\lambda\rangle,0).
  \end{align*}
  As $x$ and its translate $t^{-1}x$ are hyperspecial, $f'_\Omega = f_\Omega$ and $f'_x = f_x$ by Lemma~\ref{lm:hyperspecial}.

  Recall that $U_\Omega = \langle U_{f_\Omega,a}\rangle_{a \in \Phi^{\red}}$ (II.4.6.3). Let us show that
  $\red(U_{f_\Omega,a})$ is trivial if $\langle a,\lambda\rangle < 0$ and equals $\oU_{x,a}(\oK)$ otherwise. Note that
  $f_x^*(a) = f_x(a)+ \in \tR$ for any $a \in \Phi$, in the notation of II.4.6.9. If $\langle a,\lambda\rangle < 0$ then $f_\Omega(b) > f_x^*(b)$
  for $b \in \{a,2a\} \cap \Phi$ so that $U_{f_\Omega,a} \subset U_{f_x^*,a}$ and $\red(U_{f_\Omega,a}) = \{1\}$ as $\oGx$ is
  reductive (II.4.6.10(ii)).  Otherwise $U_{f_\Omega,a} = U_{f_x,a} = \U_{x,a}(\O)$ so that $\red(U_{f_\Omega,a}) =
  \oU_{x,a}(\oK)$.

  Putting this all together, we see that $\red(\G(\O) \cap t^{-1}\G(\O)t) = P_\lambda(\oK)$ by the rational Bruhat
  decomposition \cite[21.15]{bib:Borel} applied to $L_\lambda(\oK)$.

  To prove the final assertion, note first that $tU_{f_\Omega,a} t^{-1} = U_{f_{\Omega'},a}$, where $\Omega' = \{x, tx\}$. 
  We show that the left-hand side of~\eqref{eq:10} is contained in the right-hand side.
  It suffices to show that $t$ centralises~$\hNO$ and $U_{f_\Omega,a}$ whenever $a \in
  \Phi^{\red}$ and $\langle a,\lambda\rangle = 0$.  If $n \in \hNO$ with $\vnu(n) = w$ then $ntn^{-1} =
  n\lambda(\varpi)n^{-1} = (w\lambda)(\varpi) = \lambda(\varpi) = t$ by the above. It is a standard fact that $\im(\lambda)$
  centralises $U_a \supset U_{f_\Omega,a}$ if $\langle a,\lambda\rangle = 0$ \cite[15.4.4]{bib:Springer_LAG}.

  To prove the opposite containment in~\eqref{eq:10}, it is enough to show that the left-hand side contains $(g_+,1)$ for all $g_+ \in U_\lambda(\oK)$.
  But this is clear since we showed above that $\red(U_{f_\Omega,a}) = \oU_{x,a}(\oK)$ and $\red(U_{f_\Omega',a}) = \{1\}$ if
  $a \in \Phi^\red$ and $\langle a,\lambda\rangle > 0$.
\end{proof}

\begin{lm}[Iwahori decomposition]\label{lm:iwahori}
  Let $I \subset \G(\O)$ be the inverse image of $\oT(\oK)\oU^+(\oK)$ under the reduction map. Then $I$ is an Iwahori
  subgroup and the product map $(I\cap U^-(K)) \times (I \cap T(K)) \times (I \cap U^+(K)) \to I$ is a bijection, for any chosen order 
  of the factors. Moreover $T^-$ \(see Def.~\ref{df:t-minus}\) contracts $I \cap U^-(K)$ and expands $I \cap U^+(K)$.
\end{lm}

\begin{proof}
  By Thm.~II.4.6.33 there is a chamber $C \subset A$ with $x \in \overline C$ such that $I = \G_C^0(\O)$. Thus $I$ is an
  Iwahori subgroup.  We will use the notation of II.4.6.3. By II.4.6.7(i), $I = P^0_f = H^0 U_f$, where $f = f_C'$. Also $H^0
  = \T(\O)$ as $\T$ is connected. Note that $N_f \subset U_f \subset G(K)^1$.  Since $C$ is not contained in any walls, $N_f
  \le H$ (see the proof of Lemma~\ref{lm:fixer}). Thus $N_f \subset G(K)^1 \cap H = \T(\O)$.
  
  From I.6.4.9(iii), $I = H^0 U_f = \T(\O) U_f^+ U_f^-$. Note that $\T(\O) \subset H = \ker\nu$ normalises each $U_{a,k}$ and
  therefore $U_f^\pm$: this follows from the definitions in I.6.2. The product map is injective since $U^- \times T \times
  U^+ \to G$ is an open immersion (the big cell).
  For the final claim note that $U_f^-$ is generated by the~$U_{a,f(a)}$ for $a \in \Phi^-$ and that
  $tU_{a,k}t^{-1} \subset U_{a,k}$ for $t \in T^-$ and $a \in \Phi^-$ (II.4.2.7(2)).
\end{proof}

\begin{proof}[Proof of Lemma~\ref{lm:root_subgroup}]
  (i) Let $\Psi = \{ b \in \Phi^+: b \not\in \Z a \} \subset \Phi$. Since $a$ is simple, $\Psi$ is closed.  Thus $\prod_{b
    \in \Phi^+\nd, b \ne a} \U_{x,b} \to \U^+$ is an isomorphism (as $\O$-schemes) onto a closed subgroup scheme~$\U'$
  of~$\U^+$ (II.4.6.2). $\U'$ normal in~$\U^+$ means that the conjugation map $\U^+ \times \U' \to \U^+$ factors
  through~$\U'$, which can be checked on the generic fibre due to the $\O$-flatness of $\U^+ \times \U'$ (II.1.2.5).
  But there it is clear from $[U_b,U_c] \subset \langle U_{rb + sc} : r, s > 0 \rangle$ (condition (DR2) in I.6.1.1).
  The product map $\U_{x,a} \times \U' \to \U^+$ is an isomorphism of $\O$-schemes (II.4.6.2). As $\U'$ is normal,
  it is an isomorphism of $\O$-group schemes $\U_{x,a} \ltimes \U' \to \U^+$.

  (iii) First note that $\Gamma'_b = \Z$ for all $b \in \Phi$. This is clear when $G$ is split (II.4.2.1) and the general case
  follows either by \'etale descent (II.5.1.19) or by quasi-split descent (II.4.2.21).
  
  By II.4.5.1, $\U_{x,a}(\O) = U_{a,f_x(a)}U_{2a,f_x(2a)}$, and this equals $U_{a,f_x(a)}$ as $f_x(2a) = 2f_x(a)$ (see condition (V4) in I.6.2.1). 
  Next, from $\varphi_a(t u t^{-1}) =
  \varphi_a(u) + (\ord_K\circ a)(t)$ (II.4.2.7(2)) it follows that $tU_{a,f_x(a)}t^{-1} = U_{a,k}$, where $k =
  f_x(a)+\langle\lambda,a\rangle$. Let $l = f_x(a)+1$ so that $k$, $l \in \Gamma'_a$ and $k > l$. Then $\red(U_{a,l}) =
  \{1\}$ since $f_x^*(a) = f_x(a)+ \in \tR$, $\oGx$ is reductive, and $l > f_x(a)$ (II.4.6.10(ii)).

  Suppose first that $2a \not\in \Phi$. Then $U_a(K)$ is abelian and
  \begin{equation*}
      \sum_{U_a(K)/U_{a,k}} \psi(u_a) = \sum_{U_{a,l}\backslash U_a(K)}\, \sum_{U_{a,k}\backslash U_{a,l}} \psi(u_2u_1).
  \end{equation*}
  We claim that $U_{a,k} \subset U_{a,l}$ is a proper subgroup of $p$-power index.  This will finish the proof, since $\psi$
  is left $U_{a,l}$-invariant and the codomain~$A$ of~$\psi$ has exponent~$p$.
  Since $k$, $l \in \Gamma'_a$ and $k > l$, it follows that $U_{a,k} \subsetneq U_{a,l}$. From II.4.3.2 we see
  that $U_a(K)$ is isomorphic to the additive group of a finite (unramified) extension $L$ of~$K$. Under this isomorphism, for any $r \in \Gamma_a$,
  $U_{a,r}$ corresponds to the $\O_L$-lattice $\{x \in L: \ord_K(x) \ge r\}$. Thus the index $[U_{a,l} : U_{a,k}]$ is a power of~$p$.

  Now suppose that $2a \in \Phi$. We know that $U_{2a}(K)$ is central in~$U_a(K)$ with abelian quotient (condition (DR2) in I.6.1.1).
  Moreover, from the definitions, $U_{2a,2r} = U_{2a}(K) \cap U_{a,r}$ for all $r \in \R$.
  Note that
  \begin{equation}\label{eq:1}
    \sum_{U_a(K)/U_{a,k}} \psi(u_a) = \sum_{U_a(K)/U_{a,k}U_{2a}(K)} \psi'(u_a')
  \end{equation}
  where $\psi'(u_a') = \sum_{U_{a,k}U_{2a}(K)/U_{a,k}} \psi(u_a'u)$, which is left invariant by~$U_{a,l}$. 
  Since $U_a(K)/U_{2a}(K)$ is abelian, left and right cosets of $U_{a,k}U_{2a}(K)$ in~$U_a(K)$ coincide and we can rewrite~\eqref{eq:1} as
  \begin{align*}
    \sum_{U_{a,l}U_{2a}(K)\backslash U_a(K)}\, \sum_{U_{a,k}U_{2a}(K)\backslash U_{a,l}U_{2a}(K)} \psi'(u_2u_1)
  \end{align*}
  We claim that $U_{a,k}U_{2a}(K) \subset U_{a,l}U_{2a}(K)$ is a proper subgroup of $p$-power
  index. As in the previous case this will finish the proof.

  To see that $U_{a,k}U_{2a}(K) \subsetneq U_{a,l}U_{2a}(K)$, we show that $U_{a,k}U_{2a,2l} \subsetneq U_{a,l}$. Since $l \in \Gamma'_a$,
  we may pick $u \in U_a(K)$ such that $\varphi_a(u) = l$ and $\varphi_a(u) = \max \varphi_a(u U_{2a}(K))$. It follows that
  $u \in U_{a,l} - U_{a,k}U_{2a,2l}$.
  The index of $U_{a,k}U_{2a}(K)$ in $U_{a,l}U_{2a}(K)$ equals the index of $U_{a,k}/U_{2a,2k}$ in $U_{a,l}/U_{2a,2l}$. The group
  $U_a(K)/U_{2a}(K)$ is isomorphic to the additive group of a finite-dimensional $K$-vector space and for any $r \in \Gamma'_a$,
  $U_{a,r}/U_{2a,2r}$ corresponds to an $\O$-lattice under this isomorphism (II.4.3.7, II.4.3.5 with $k = r$, $l = 2r \in \Gamma'_{2a}$).
  Thus the index of $U_{a,k}U_{2a}(K)$ in $U_{a,l}U_{2a}(K)$ is a $p$-power.
\end{proof}

\bibliography{satake_mod_p}
\bibliographystyle{halpha}

\end{document}